\newcommand{\dsp}{\renewcommand{\baselinestretch}{1.2}}

\documentclass[pdflatex,fleqn]{article}
\usepackage{}
\usepackage{xcolor}
\usepackage{mathrsfs}
\usepackage{amsfonts}
\usepackage{amssymb}
\usepackage{latexsym}
\usepackage{amsmath}
\usepackage{amscd}
\usepackage{graphicx}
\usepackage{hyperref}
\usepackage{hypernat}
\usepackage{cite}

\usepackage{geometry,amsthm,graphics,amssymb,amsmath,enumerate,latexsym,tabularx,shapepar}
\usepackage[all,2cell,dvips]{xy}\UseAllTwocells\SilentMatrices

\usepackage{extarrows}
\date{}

\dsp

\title{\bf A reduction theorem for $AH$ algebras with the ideal property}

\author{
Guihua Gong, Chunlan Jiang,  Liangqing Li and Cornel Pasnicu
}

\begin{document}

\large

\maketitle

\noindent{\underline {\noindent{ }\hspace {157mm}}}

\noindent{\bf Abstract}

Let $A$ be an $AH$ algebra, that is, $A$ is the inductive limit $C^{*}$-algebra of $$A_{1}\xrightarrow{\phi_{1,2}}A_{2}\xrightarrow{\phi_{2,3}}A_{3}\longrightarrow\cdots\longrightarrow A_{n}\longrightarrow\cdots$$
with $A_{n}=\bigoplus_{i=1}^{t_{n}}P_{n,i}M_{[n,i]}(C(X_{n,i}))P_{n,i}$, where $X_{n,i}$ are compact metric spaces, $t_{n}$ and $[n,i]$ are positive integers, and $P_{n,i}\in M_{[n,i]}(C(X_{n,i}))$ are projections. Suppose that $A$ has the ideal property: each closed two-sided ideal of $A$ is generated by the projections inside the ideal, as a closed two-sided  ideal. Suppose that $\sup_{n,i}dim(X_{n,i})<+\infty$. In this article, we prove that $A$ can be written as the inductive limit of $$B_{1}\longrightarrow B_{2}\longrightarrow\cdots\longrightarrow B_{n}\longrightarrow\cdots,$$
where $B_{n}=\bigoplus_{i=1}^{s_{n}}Q_{n,i}M_{\{n,i\}}(C(Y_{n,i}))Q_{n,i}$, where $Y_{n,i}$ are $\{pt\}, [0,1], S^{1}, T_{\uppercase\expandafter{\romannumeral2}, k}, T_{\uppercase\expandafter{\romannumeral3}, k}$ and $S^{2}$ (all of them are connected simplicial complexes of dimension at most three), $s_{n}$ and $\{n,i\}$ are positive integers and $Q_{n,i}\in M_{\{n,i\}}(C(Y_{n,i}))$ are projections. This theorem unifies and generalizes the reduction theorem for real rank zero $AH$ algebras due to Dadarlat and Gong ([D], [G3] and [DG]) and the reduction theorem for simple $AH$ algebras due to Gong (see [G4]).

\vspace{5mm}

 
 \newpage
 
\noindent\textbf{\S1. Introduction}

\vspace{3mm}
Successful classification results have been obtained for real rank zero $AH$ algebras (see [Ell1], [Lin1-3] [EG1-2], [EGLP], [D], [G1-3], [DG]) and simple $AH$ algebras (see [Ell2-3], [Li1-3], [G4], [EGL1-2]) in the case of no dimension growth (this condition can be relaxed to a certain slow dimension growth condition). To unify these two classification theorems, we will consider $AH$ algebras with the ideal property (see [Ji-Jiang] and [GJLP]). This article is a continuation of the paper [GJLP]---we obtain the reduction theorem for arbitrary $AH$ algebras $A$ (of no dimension growth) with the ideal property. That is, we remove the restriction that  $K_{*}(A)$  is torsion free in the paper [GJLP]. Since we do not assume that $K_{*}(A)$ is torsion free,  we must involve higher dimensional spaces such as $T_{\uppercase\expandafter{\romannumeral2}, k}$, 
$T_{\uppercase\expandafter{\romannumeral3}, k}$, and $S^2$ in our reduction theorem (see below). This makes the main result of this paper  much more difficult to prove than the one in [GJLP], as $C(X)$ is not stably generated when $\dim(X)\geq 2$. 

The ideal property is a property of structural interest for a $C^*$-algebra. Many interesting and important $C^*$-algebras have the ideal property. It was proved by Cuntz-Echterhoff-Li that semigroup $C^*$-algebras of $ax + b$-semigroups over Dedekind domains have the ideal property ([CEL]). A generalization of this result appeared in [L]. Interesting examples of crossed product $C^*$-algebras with the ideal property could be found, e.g., in [Pa-Ph1] and [Pa-Ph2]. Other important results involving the ideal property have been proved, e.g., in [Pa-R1-2], [Pa-Ph1] and [Pa-Ph2]. Many $C^*$-algebras coming from $\mathbb{Z}$ dynamical systems on compact metric spaces are $AH$ algebras (see [EE],  [Phi1-2], [Lin4] and [LinP]).

An $AH$ algebra is a nuclear $C^{*}$-algebra of the form $A=\lim\limits_{\longrightarrow}(A_{n}, \phi_{n,m})$ with $A_{n}=\bigoplus_{i=1}^{t_{n}}P_{n,i}M_{[n,i]}C(X_{n,i})P_{n,i}$, where $X_{n,i}$ are compact metric spaces, $t_{n}$,  $[n,i]$ are positive integers, $M_{[n,i]}(C(X_{n,i}))$ are algebras of $[n,i]\times[n,i]$ matrices with entries in $C(X_{n,i})$---the algebra of complex - valued continuos functions on $X_{n,i}$---, and finally $P_{n,i}\in M_{[n,i]}(C(X_{n,i}))$ are projections (see [Bla]). Let $T_{\uppercase\expandafter{\romannumeral2}, k}$ (and $T_{\uppercase\expandafter{\romannumeral3}, k}$) be a connected finite simplicial complex with $H^{1}(T_{\uppercase\expandafter{\romannumeral2}, k})=0$ and $H^{2}(T_{\uppercase\expandafter{\romannumeral2}, k})=\mathbb{Z}/k\mathbb{Z}$ (and $H^{1}(T_{\uppercase\expandafter{\romannumeral3}, k})=0=H^{2}(T_{\uppercase\expandafter{\romannumeral3}, k})$ and $H^{3}(T_{\uppercase\expandafter{\romannumeral3}, k})=\mathbb{Z}/k\mathbb{Z}$, respectively). Recall that the unit circle is denoted by $S^{1}$ and  the   2-dimensional unit sphere is denoted by $S^{2}$. In this article, we will prove that an $AH$ algebra with ideal property and no dimension growth can be rewritten as an $AH$ inductive limit with the spaces $X_{n,i}$ being $\{pt\}, [0,1], S^{1}, T_{\uppercase\expandafter{\romannumeral2}, k}$, $T_{\uppercase\expandafter{\romannumeral3}, k}$ and $S^{2}$. (For the background information, we refer the readers to [GJLP].) The result in this paper plays an essential  role in the classification of the $AH$ algebras with the ideal property and with no dimension growth (see [GJL]).


\vspace{3mm}

\noindent\textbf{\S2. Notation and prelimilary}

\vspace{-2.4mm}

The following notations are quoted from [G4] or [GJLP].

\noindent\textbf{2.1.} If $A$ and $B$ are two $C^{\ast}$-algebras, we use $Map(A,B)$ to denote \textbf{the space of all linear, completely positive $\ast$-contractions} from $A$ to $B$. If both $A$ and $B$ are unital, then $Map(A,B)_{1}$ will denote the subset of $Map(A,B)$ consisting of unital maps. By word "map", we shall mean a linear, completely positive $\ast$-contraction between $C^{\ast}$-algebras, or else we shall mean a continuous map between topological spaces, which one will be clear from the context.

\textbf{By a homomorphism between $C^{\ast}$-algebras, will be meant a $\ast$-homomorphism.} Let $Hom(A,B)$ denote the \textbf{space of all the homomorphisms} from $A$ to $B$. Similarly, if both $A$ and $B$ are unital, let $Hom(A,B)_{1}$ denote the subset of $Hom(A,B)$ consisting of all the untial homomorphisms.

\noindent\textbf{Definition 2.2.} Let $G\subset A$ be a finite set and $\delta>0$. We  shall say that $\phi\in Map(A,B)$ is $G-\delta$ \textbf{multiplicative} if
$$\parallel\phi(ab)-\phi(a)\phi(b)\parallel<\delta$$
for all $a,b\in G$.

\noindent\textbf{2.3.} In the notation for an inductive system $(A_{n},\phi_{n,m})$, we understand that $\phi_{n,m}=\phi_{m-1,m}\circ\phi_{m-2,m-1}\circ\cdots\circ\phi_{n,n+1},$
where all $\phi_{n,m}:A_{n}\rightarrow A_{m}$ are homomorphisms.

We shall assume that, for any summand $A^{i} _{n}$ in the direct sum $A_{n}=\bigoplus^{t_{n}}_{i=1}A^{i} _{n}$, necessarily
$\phi_{n,n+1}(\textbf{1}_{A^{i} _{n}})\neq0$, since, otherwise, we could simply delete $A^{i} _{n}$ from $A_{n}$ without changing the limit algebra.

\noindent\textbf{2.4.} If $A_{n}=\bigoplus_{i}A^{i} _{n}$ and $A_{m}=\bigoplus_{j}A^{j} _{m}$, we use $\phi^{i,j}_{n,m}$ to denote the partial map of $\phi_{n,m}$ from the
$i$-th block $A^{i} _{n}$ of $A_{n}$ to the $j$-th block $A^{j} _{m}$ of $A_{m}$.

\noindent\textbf{2.5.} By 2.3 of [Bla] and Theorem 2.1 of [EGL2], we know that any $AH$ algebra can be  written as an inductive limit $A=lim(A_{n}=\bigoplus_{i=1}^{t_{n}}P_{n,i}M_{[n,i]}C(X_{n,i})P_{n,i}, \phi_{n,m})$, where $X_{n,i}$ are finite simplicial complexes and $\phi_{n,m}$ are injective. In this article we will assume that $X_{n,i}$ are (path) connected finite simplicial complexes and $\phi_{n,m}$ are injective.

It is well known that, for any connected finite simplicial complex $X$, there is a metric $d$ on $X$ with the following property: for any $x\in X$ and $\eta>0$, the $\eta$-ball centered in $x$, $B_{\eta}(x)=\{x^{\prime}\in X|~ d(x^{\prime},x)<\eta\}$ is path connected. So in this article, we will always assume that the metric on a connected simplicial complex has this property.

\noindent\textbf{2.6.} In this article, we assume that the  inductive limit system satisfies the no dimension growth condition: there is an $M\in \mathbb{N}$ such that for any $n,i,$
$$dim(X_{n,i})\leqslant M.$$

The condition could be relaxed to a so called very slow dimension growth for our main theorem. Since the proof of this slightly more general case is quite tedious, we will leave it to a subsequent paper.

\noindent\textbf{2.7.} By 1.3.3 of [G4], any $AH$ algebra $$A=lim(A_{n}=\bigoplus_{i=1}^{t_{n}}P_{n,i}M_{[n,i]}(C(X_{n,i}))P_{n,i}, \phi_{n,m})$$
is isomorphic to a limit corner subalgebra (see Definition 1.3.2 of [G4] for this concept) of $\tilde{A}=\lim\limits_{\longrightarrow}(\tilde{A}_{n}=\bigoplus_{i=1}^{t_{n}}M_{[n,i]}(C(X_{n,i})), \tilde{\phi}_{n,m})$---an inductive limit of full matrix algebras over $X_{n,i}$. Once we prove that $\tilde{A}$ is an inductive limit of homogeneous algebras over the spaces in  the following list: $\{pt\}, [0,1], S^{1}, T_{\uppercase\expandafter{\romannumeral2}, k}, T_{\uppercase\expandafter{\romannumeral3}, k}, S^{2}$, then $A$ itself is also an inductive limit of such kind. Therefore in this article, we will assume $A$ itself is an inductive limit of finite direct sum of full matrix algebras over $X_{n,i}$. That is $P_{n,i}=\textbf{1}_{M_{[n,i]}(C(X_{n,i}))}$.

\noindent\textbf{2.8.} Let $Y$ be a compact metrizable space. Let $P\in M_{k_{1}}(C(Y))$ be a projection with rank($P$)=$k\leqslant k_{1}$. For each $y$, there is a unitary $u_{y}\in M_{k_{1}}(\mathbb{C})$ (depending on $y$) such that
$$P(y)=u_{y}\left(
            \begin{array}{cccccc}
              1 & \ & \ & \ & \ &\ \\
              \ & \ddots & \ & \ & \ &\ \\
              \ & \ & 1 & \ & \ &\ \\
              \ & \ & \ & 0 & \ &\ \\
              \ & \ & \ & \ & \ddots \ & \ \\
               \ & \ & \ & \ & \ & 0 \  \\
                \end{array}
              \right)u^{*}_{y},$$
where there are $k$ $1^{'}s$ on the diagonal. If the unitary $u_{y}$ can be chosen to be continuous in $y$, then $P$ is called a \textbf{trivial projection}.

It is well known that any projection $P\in M_{k_{1}}(C(Y))$ is locally trivial. That is, for any $y_{0}\in Y$, there is an open set $U_{y_{0}}\ni y_{0}$, and there is a continuous unitary-valued function
$$u: U_{y_{0}}\longrightarrow M_{k_{1}}(\mathbb{C})$$
such that the above equation holds for $u(y)$ (in place of $u_{y}$) for any $y\in U_{y_{0}}$.

If $P$ is trivial, then $P M_{k_{1}}(C(X))P\cong M_{k}(C(X))$.

\noindent\textbf{2.9.} Let $X$ be a compact metrizable space and $\psi: C(X)\longrightarrow P M_{k_{1}}(C(Y))P$ be a unital homomorphism. For any given point $y\in Y$, there are points
$$x_{1}(y),x_{2}(y),\cdots,x_{k}(y)\in X,$$
and a unitary $U_{y}\in M_{k_{1}}(\mathbb{C})$ such that
$$\psi(f)(y)=P(y)U_{y}\left(
            \begin{array}{cccccc}
              f(x_{1}(y)) & \ & \ & \ & \ &\ \\
              \ & \ddots & \ & \ & \ &\ \\
              \ & \ & f(x_{k}(y)) & \ & \ &\ \\
              \ & \ & \ & 0 & \ &\ \\
              \ & \ & \ & \ & \ddots \ & \ \\
               \ & \ & \ & \ & \ & 0 \  \\
                \end{array}
              \right)U^{*}_{y}P(y)\in P(y)M_{k_{1}}(\mathbb{C})P(y),$$
for all $f\in C(X)$. Equivalently, there are $k$ rank one orthogonal projections $p_{1},p_{2},\cdots,p_{k}$ with $\sum^{k}_{i=1}p_{i}(y)=P(y)$ and $x_{1}(y),x_{2}(y),\cdots,x_{k}(y)\in X$, such that
$$\psi(f)(y)=\sum\limits^{k}\limits_{i=1}f(x_{i}(y))p_{i}(y),~~~~\forall f\in C(X).$$

Let us denote the set $\{x_{1}(y),x_{2}(y),\cdots,x_{k}(y)\}$, counting multiplicities, by SP$\psi_{y}$ (see [Pa1]). In other words, if a point is repeated in the diagonal of the above matrix, it is included with the same multiplicity in SP$\psi_{y}$.\textbf{ We shall call SP$\psi_{y}$ the spectrum of $\psi$ at the point} $y$. Let us define the \textbf{spectrum} of $\psi$, denoted by SP$\psi$, to be the closed subset
$$SP\psi:= \overline{\bigcup\limits_{y\in Y}SP\psi_{y}}(=\bigcup\limits_{y\in Y}SP\psi_{y})\subset X.$$
Alternatively, SP$\psi$ is the complement of the specturm of the kernel of $\psi$, considered as a closed ideal of $C(X)$. The map $\psi$ can be factored as
$$C(X)\xrightarrow{i^{\ast}}C(SP\psi)\xrightarrow{\psi_{1}}PM_{k_{1}}(C(Y))P$$
with $\psi_{1}$ an injective homomorphism, where $i$ denotes the inclusion $SP\psi\hookrightarrow X.$

Also, if $A=PM_{k_{1}}(C(Y))P$, then we shall call the space $Y$ the spectrum of the algebra $A$, and write $SP A=Y(=SP(id))$.

\noindent\textbf{2.10.} In 2.9, if we group together all the repeated points in $\{x_{1}(y),x_{2}(y),\cdots,x_{k}(y)\}$, and sum their corresponding projections, we can write
$$\psi(f)(y)=\sum\limits^{l}\limits_{i=1}f(\lambda_{i}(y))P_{i},~~~~(l\leqslant k),$$
where $\{\lambda_{1}(y),\lambda_{2}(y),\cdots,\lambda_{l}(y)\}$ is equal to $\{x_{1}(y),x_{2}(y),\cdots,x_{k}(y)\}$ as a set, but $\lambda_{i}(y)\neq \lambda_{j}(y)$ if $i\neq j$; and each $P_{i}$ is the sum of the projections corresponding to $\lambda_{i}(y)$. If $\lambda_{i}(y)$ has multiplicity $m$ (i.e., it appears $m$ times in $\{x_{1}(y),x_{2}(y),\cdots,x_{k}(y)\}$), then rank($P_{i}$)=$m$.

\noindent\textbf{2.11.} Set $P^{k}(X)=\underbrace{X\times X\times\cdots\times X}\limits_{k}/\thicksim$, where the equivalence relation $\thicksim$ is defined by $(x_{1},x_{2},\cdots,x_{k})\thicksim(x^{\prime}_{1},x^{\prime}_{2},\cdots,x^{\prime}_{k})$ if there is a permutation $\sigma$ of $\{1,2,\cdots,k\}$ such that
$x_{i}=x^{\prime}_{\sigma(i)}$, for each $1\leq i\leq k$. A metric $d$ on $X$ can be extended to a metric on $P^{k}(X)$ by $$d([x_{1},x_{2},\cdots,x_{k}],[x^{\prime}_{1},x^{\prime}_{2},\cdots,x^{\prime}_{k}])=\min\limits_{\sigma}\max\limits_{1\leq i\leq k}d(x_{i},x^{\prime}_{\sigma(i)}),$$
where $\sigma$ is taken from the set of all permutations, and $[x_{1},x_{2},\cdots,x_{k}]$ denotes the equivalence class in $P^{k}(X)$ of $(x_{1},x_{2},\cdots,x_{k})$. Two k-tuples of (possible repeating) points $\{x_{1},x_{2},\cdots,x_{k}\}\subset X$ and $\{x^{\prime}_{1},x^{\prime}_{2},\cdots,x^{\prime}_{k}\}\subset X$ are said to be paired within $\eta$ if
$$d([x_{1},x_{2},\cdots,x_{k}],[x^{\prime}_{1},x^{\prime}_{2},\cdots,x^{\prime}_{k}])<\eta$$ when one regards $(x_{1},x_{2},\cdots,x_{k})$ and $(x^{\prime}_{1},x^{\prime}_{2},\cdots,x^{\prime}_{k})$ as two points in $P^{k}(X)$.

\noindent\textbf{2.12.} Let $\psi: C(X)\longrightarrow P M_{k_{1}}(C(Y))P$ be a unital homomorphism as in 2.9. Then
$$\psi^{\ast}: y\mapsto SP \psi_{y}$$
defines a map $Y\longrightarrow P^{k}(X)$, if one regards $SP \psi_{y}$ as an element of $P^{k}(X)$. This map is continuous. In terms of this map and the metric $d$, let us define the \textbf{spectral variation} of $\psi$:\\
$$SPV(\psi):= the\; diamerter\; of\; image\; of\; \psi^{\ast}.$$

\noindent\textbf{Definition 2.13}~~We shall call $P_{i}$ in 2.10 the spectral projection of $\phi$ at $y$ with respect to the spectral element $\lambda_{i}(y)$. For a subset $X_{1}\subset X$, we shall call
$$\sum\limits_{\lambda_{i}(y)\in X_{1}}P_{i}$$
the spectral projection of $\phi$ at $y$ corresponding to the subset $X_{1}$ (or with respect to the subset $X_{1}$).\\

\noindent\textbf{2.14.} Let $\phi: M_{k}(C(X))\longrightarrow P M_{l}(C(Y))P$ be a unital homomorphism. Set $\phi(e_{11})=p$, where $e_{11}$ is the canonical matrix unit corresponding to the upper left corner. Set
$$\phi_{1}=\phi\mid_{e_{11}M_{k}(C(X))e_{11}}: C(X)\longrightarrow p M_{l}(C(Y))p.$$
Then $P M_{l}(C(Y))P$ can be identified with $p M_{l}(C(Y))p\otimes M_{k}$ in such a way that
$$\phi=\phi_{1}\otimes id_{k}.$$
Let us define
$$SP \phi_{y}:=SP(\phi_{1})_{y},$$
$$SP \phi:=SP\phi_{1},$$
$$SPV(\phi):=SPV(\phi_{1}).$$

Suppose that $X$ and $Y$ are connected. Let $Q$ be a projection in $M_{k}(C(X))$ and $\phi:QM_{k}(C(X))Q\longrightarrow PM_{l}(C(Y))P$ be a unital map. By the Dilation lemma(Lemma 2.13 of [EG2]), there are an $n$, a projection $P_{1}\in M_{n}(C(Y))$, and a unital homomorphism
$$\widetilde{\phi}:M_{k}(C(X))\longrightarrow P_{1}M_{n}(C(Y))P_{1}$$
such that
$$\phi=\widetilde{\phi}\mid_{QM_{k}(C(X))Q}.$$
(Note that this implies that $P$ is a subprojection of $P_{1}.$) We define:
$$SP \phi_{y}:=SP\widetilde{\phi}_{y},$$
$$SP \phi:=SP\widetilde{\phi},$$
$$SPV(\phi):=SPV(\widetilde{\phi}).$$
(Note that these definitions do not depend on the choice of the dilation $\widetilde{\phi}$.)

\noindent\textbf{2.15.} Let $\phi: M_{k}(C(X))\longrightarrow P M_{l}(C(Y))P$ be a (not necessarily unital) homomorphism, where $X$ and $Y$ are connected finite simplicial complexes. Then
$$\#(SP \phi_{y})=\frac{rank\phi(\textbf{1}_{k})}{rank(\textbf{1}_{k})},~~~for\; any\; y\in Y,$$
where again $\#(\cdot)$ denotes the number of elements in the set counting multiplicity. It is also true that for any nonzero projection $p\in M_{k}(C(X))$, $\#(SP \phi_{y})=\frac{rank\phi(p)}{rank(p)}$.

\noindent\textbf{2.16.} Let $X$ be a compact connected space and let $Q$ be a projection of rank $n$ in $M_{N}(C(X))$. The \textbf{weak variation of a finite set} $F\subset QM_{N}(C(X))Q$ is defined by
$$\omega(F)=\sup\limits_{\Pi_{1},\Pi_{2}}\inf\limits_{u\in U(n)}\max\limits_{a\in F}\parallel u\Pi_{1}(a)u^{\ast}-\Pi_{2}(a)\parallel,$$
where $\Pi_{1},\Pi_{2}$ run through the set of irreducible representations of $QM_{N}(C(X))Q$ into $M_{n}(\mathbb{C})$.

Let $X_{i}$ be compact connected spaces and $Q_{i}\in M_{n_{i}}(C(X_{i}))$ be projections. For a finite set $F\subset\bigoplus_{i}Q_{i}M_{n_i}(C(X_{i}))Q_{i}$, define the \textbf{weak variation} $\omega(F)$ to be $\max_{i} \omega(\pi_{i}(F))$, where $\pi_{i}: \bigoplus_{j}Q_{j}M_{n_{j}}(C(X_{j}))Q_{j}\longrightarrow Q_{i}M_{n_{i}}(C(X_{i}))Q_{i}$ is the natural projection map onto the $i$-th block.

The set $F$ is said to be \textbf{weakly approximately constant to within} $\varepsilon$ if $\omega(F)<\varepsilon$.

\noindent\textbf{2.17.} The following notations will be frequently used in this article.\\
(a) As in 2.15, we use notation $\#(\cdot)$ to denote the cardinal number of the set counting multiplicity.\\
(b) For any metric space $X$, any $x_{0}\in X$ and any $c>0$, let
$$B_{c}(x_{0}):= \{x\in X \mid d(x,x_{0})<c\}$$
denote the open ball with radius $c$ and center $x_{0}$.\\
(c) Suppose that $A$ is a $C^{\ast}$-algebra, $B\subset A$ is a sub-$C^{*}$-algebra, $F\subset A$ is a (finite) subset and let $\varepsilon>0$. If for each element $f\in F$, there is an element $g\in B$ such that $\parallel f-g\parallel<\varepsilon$, then we shall say that $F$ is approximately contained in $B$ to within $\varepsilon$, and denote this by $F\subset_{\varepsilon} B$.\\
(d) Let $X$ be a compact metric space. For any $\delta>0$, a finite set $\{x_{1},x_{2},\cdots,x_{n}\}$ is said to be $\delta$-dense in $X$, if for any $x\in X$, there is $x_{i}$ such that dist$(x,x_{i})<\delta$.\\
(e) We shall use $\bullet$ to denote any possible positive integer. To save notation, $a_{1},a_{2},\cdots$ may be used for a finite sequence if we do not care how many terms are in the sequence. Similarly, $A_{1}\cup A_{2}\cup\cdots$ or $A_{1}\cap A_{2}\cap\cdots$ may be used for a finite union or a finite intersection. If there is a danger of confusion with an infinite sequence, union or intersection, we will write them as $a_{1},a_{2},\cdots,a_{\bullet}$, $A_{1}\cup A_{2}\cup\cdots\cup A_{\bullet}$, $A_{1}\cap A_{2}\cap\cdots\cap A_{\bullet}$.\\
(f) In this paper, we often use $\textbf{1}$ to denote the units of different unital $C^{\ast}$-algebras. In particular, if $\textbf{1}$ appears in $\phi(\textbf{1})$, where $\phi$ is a homomorphism, then $\textbf{1}$ is the unit of the domain algebra. For example for a homomorphism $\phi: \bigoplus^{r}_{i=1}A^{i}\longrightarrow B$, then $\textbf{1}$ in $\phi(\textbf{1})$ means $\textbf{1}_{\bigoplus^{r}_{i=1}A^{i}}$ and $\textbf{1}$ in $\phi^{i}(\textbf{1})$ means $\textbf{1}_{A^{i}}$.\\
(g) For any two projections $p,q\in A$, we use the notation $[p]\leqslant[q]$ to denote that $p$ is unitarily equivalent to a subprojection of $q$. And we use $p\thicksim q$ to denote that $p$ is unitarily equivalent to $q$.

\noindent\textbf{2.18.} For any $\eta>0$, $\delta>0$, a unital homomorphism $\phi: PM_{k}(C(X))P\longrightarrow Q M_{k^{'}}(C(Y))Q$ is said to have the property $sdp(\eta,\delta)$(spectral distribution property with respect to $\eta$ and $\delta$) if for any $\eta$-ball $B_{\eta}(x)$ and any point $y\in Y$,
$$\#(SP \phi_{y}\cap B_{\eta}(x))\geqslant\delta\#(SP \phi_{y})(=\delta\frac{rank(Q)}{rank(P)})$$
counting multiplicity. Any homomorphism
$$\phi:\bigoplus_{i}P_{i}M_{k_{i}}(C(X_{i}))P_{i}\longrightarrow \bigoplus_{j}Q_{j}M_{l_{j}}(C(Y_{j}))Q_{j}$$
is said to have the property $sdp(\eta,\delta)$ if each partial map
$$\phi^{i,j}:P_{i}M_{k_{i}}(C(X_{i}))P_{i}\longrightarrow \phi^{i,j}(P_{i})M_{l_{j}}(C(Y_{j}))\phi^{i,j}(P_{i})$$
has the property $sdp(\eta,\delta)$ as a unital homomorphism. Note that by definition, a nonunital homomorphism $\phi: M_{k}(C(X))\longrightarrow M_{l}(C(Y))$ has the property $sdp(\eta,\delta)$  if the corresponding unital map
$$\phi: M_{k}(C(X))\longrightarrow \phi(\textbf{1}_{k})M_{l}(C(Y))\phi(\textbf{1}_{k})$$
has the property $sdp(\eta,\delta)$.

The following Lemma is Lemma 2.8 of [GJLP].

\noindent\textbf{Lemma 2.19.} Let $A=\lim\limits_{\longrightarrow}(A_{n}=\bigoplus_{i=1}^{t_{n}}M_{[n,i]}(C(X_{n,i})), \phi_{n,m})$ be an $AH$ algebra with the ideal property. For $A_{n}$ and any $\eta>0$, there exist a $\delta>0$, a positive integer $m>n$, connected finite simplicial complexes $Z^{1}_{i},Z^{2}_{i},\cdots,Z^{\bullet}_{i}\subset X_{n,i}, i=1,2,\cdots,t_{n}$, and a homomorphism
$$\phi: B=\bigoplus\limits_{i=1}\limits^{t_{n}}\bigoplus\limits_{s}M_{[n,i]}(C(Z^{s}_{i}))\longrightarrow A_{m}$$
such that\\
(1) $\phi_{n,m}$ factors as $A_{n}\xrightarrow{\pi}B\xrightarrow{\phi}A_{m}$, where $\pi$ is defined by
$$\pi(f)=(f\mid_{Z^{1}_{i}},f\mid_{Z^{2}_{i}},\cdot\cdot\cdot,f\mid_{Z^{\bullet}_{i}})\in \bigoplus\limits_{s}M_{[n,i]}(C(Z^{s}_{i}))\subset B,$$
for any $f\in M_{[n,i]}(C(X_{n,i})).$\\
(2) The homomorphism $\phi$ satisfies the dichotomy condition ($\ast$): for each $Z^{s}_{i}$, the partial map
$$\phi^{(i,s),j}: M_{[n,i]}(C(Z^{s}_{i}))\longrightarrow A_{m}^j$$
either has the property $sdp(\frac{\eta}{32},\delta)$ or is the zero map. Furthermore for any $m^{'}>m$, each partial map of $\phi_{m,m^{'}}\circ\phi$ satisfies the dichotomy condition $(\ast)$: either it has the property $sdp(\frac{\eta}{32},\delta)$ or is the zero map.

The following result is quoted from [Pa2].

\noindent\textbf{Lemma 2.20.} Let $A=lim(A_{n}=\bigoplus_{i=1}^{t_{n}}M_{[n,i]}(C(X_{n,i})), \phi_{n,m})$ be an $AH$ inductive limit with the ideal property and with no dimension growth (as in 2.5, we assume $X_{n,i}$ path connected). For any $A_{n}$, finite  set $F_{n}=\bigoplus F^{i}_{n}\subset A_{n}$, $\varepsilon>0$, and positive integer $L$, there is an $A_{m}$, such that for each pair $(i,j)$, one of  the following conditions holds \\
(i) $\frac{rank(\phi^{i,j}_{n,m}(\textbf{1}_{A^{i}_{n}}))}{rank(\textbf{1}_{A^{i}_{n}})}\geqslant L$, or\\
(ii) there is a homomorphism
$$\psi: A^{i}_{n}\longrightarrow \phi^{i,j}_{n,m}(\textbf{1}_{A^{i}_{n}})A^{j}_{m}\phi^{i,j}_{n,m}(\textbf{1}_{A^{i}_{n}})$$
with finite dimensional image such that $\psi$ is homotopic to $\phi^{i,j}_{n,m}$ and
$$\parallel \psi(f)-\phi^{i,j}_{n,m}(f)\parallel<\varepsilon~~~~~\forall f\in F^{i}_{n}.$$
The following proposition is Theorem 2.12 of [GJLP].

\noindent\textbf{Proposition 2.21.}
Let $A=\lim\limits_{\longrightarrow}(A_{n}=\bigoplus_{i=1}^{t_{n}}M_{[n,i]}(C(X_{n,i})), \phi_{n,m})$
be an $AH$ algebra with the ideal property and with no dimension growth. For any $A_{n}$, finite set $F=\bigoplus_{i=1}^{t_{n}}F^{i}_{n}\subset A_{n}$, positive integer $J$ and $\varepsilon>0$, there exists $m$ and there exist projections $Q_{0},Q_{1},Q_{2}\in A_{m}$ with $Q_{0}+Q_{1}+Q_{2}=\phi_{n,m}(\textbf{1}_{A_{n}})$, a unital map $\psi_{0}\in Map(A_{n},Q_{0}A_{m}Q_{0})_{1}$ and two unital homomorphisms $\psi_{1}\in Hom(A_{n},Q_{1}A_{m}Q_{1})_{1}$, $\psi_{2}\in Hom(A_{n},Q_{2}A_{m}Q_{2})_{1}$ such that the following statement are true \\
(1) $\parallel \phi_{n,m}(f)-(\psi_{0}(f)\oplus\psi_{1}(f)\oplus\psi_{2}(f))\parallel<\varepsilon$, for all $f\in F$\\
(2) $\omega((\psi_{0}\oplus\psi_{1})(F))<\varepsilon$\\
(3) The homomorphism $\psi_{2}$ factors through $C$---a finite direct sum of matrix algebras over $C[0,1]$, or $\mathbb{C}$ as
$$\psi_{2}: A_{n}\xrightarrow{\xi_{1}}C\xrightarrow{\xi_{2}}Q_{2}A_{m}Q_{2}$$
where $\xi_{1}, \xi_{2}$ are unital homomorphisms\\
(4) Let $\psi^{i,j}_{0}:
A^{i}_{n}\longrightarrow\psi^{i,j}_{0}(\textbf{1}_{A^{i}_{n}})A^{j}_{m}\psi^{i,j}_{0}(\textbf{1}_{A^{i}_{n}})$ and
$\psi^{i,j}_{1}: A^{i}_{n}\longrightarrow\psi^{i,j}_{1}(\textbf{1}_{A^{i}_{n}})A^{j}_{m}\psi^{i,j}_{1}(\textbf{1}_{A^{i}_{n}})$ be the corresponding partial maps of $\psi_{0}$ and $\psi_{1}$. For each pair $(i,j)$, one of the following is true.\\
(i) Both $\psi^{i,j}_{0}$ and $\psi^{i,j}_{1}$ are zero, or\\
(ii) $\psi^{i,j}_{1}$ is a homomorphism with finite dimensional image and for each non zero projections $e\in A^{i}_{n}$ (including any rank 1 projection)$$[\psi^{i,j}_{1}(e)]>J[\psi^{i,j}_{0}(\textbf{1}_{A^{i}_{n}})](\in K_{0}(A_{m}^j)).$$

Furthermore, we can assume that $Q^{j}_{0}, Q^{j}_{1}$ are trivial projections in $A^{j}_{m}$.

The following Corollary is also in [GJLP],

\noindent\textbf{Corollary 2.22.} We use the notation from 2.21. For any $A_{n}$, any projection $P=\bigoplus P^{i}\in\bigoplus A^{i}_{n}$, any finite set $F=\bigoplus F^{i}\in\bigoplus P^{i}A^{i}_{n}P^{i}=PA_{n}P$, any positive integer $J$, and any number $\varepsilon>0$, there are an $A_{m}$, mutually orthogonal projections $Q_{0},Q_{1},Q_{2}\in A_{m}$ with $Q_{0}+Q_{1}+Q_{2}=\phi_{n,m}(\textbf{1}_{A_{n}})$, a unital map $\psi_{0}\in Map(A_{n},Q_{0}A_{m}Q_{0})_{1}$ and two unital homomorphisms $\psi_{1}\in Hom(A_{n},Q_{1}A_{m}Q_{1})_{1}$, $\psi_{2}\in Hom(A_{n},Q_{2}A_{m}Q_{2})_{1}$ such that
for each pair $(i,j)$, $\psi^{i,j}_{0}(P^{i})$ and $\psi^{i,j}_{0}(\textbf{1}_{A_{n}^{i}}-P^{i})$ are mutually orthogonal projections and there is an approximate decomposition of $\phi{'}_{n,m}:=\phi_{n,m}\mid_{PA_{n}P}$ as a direct sum of $\psi{'}_{0}:=\psi_{0}\mid_{PA_{n}P},\psi{'}_{1}:=\psi_{1}\mid_{PA_{n}P}$ and $\psi{'}_{2}:=\psi_{2}\mid_{PA_{n}P}$, satisfying the following conditions:\\
(1) $\parallel \phi{'}_{n,m}(f)-(\psi{'}_{0}(f)\oplus\psi{'}_{1}(f)\oplus\psi{'}_{2}(f))\parallel<\varepsilon$, for all $f\in F$.\\
(2) $\psi{'}_{1}$ has finite dimensional image and $\psi{'}_{2}$ factors through a finite direct sum of matrix algebras over $C[0,1]$ or $\mathbb{C}$.\\
(3) If $\psi'^{i,j}_{0}\neq 0$, then for any nonzero projection $e\in P^{i}A^{i}_{n}P^{i}, [\psi'^{i,j}_{1}(e)]>J[\psi'^{i,j}_{0}(P^{i})](\in K_0( A^{j}_{m}))$.\\
(4) $\psi{'}_{0}$ is $F-\varepsilon$ multiplicative.

\noindent\textbf{2.23.} Let $X$ be a connected finite simplicial complex, $A=M_{k}(C(X))$. A unital $\ast$-monomorphism $\phi:A\longrightarrow M_{l}(A)$ is called a $\textbf{(unital) simple embedding} $ if it is homotopic to the homomorphism $id\oplus\lambda$, where $\lambda: A\longrightarrow M_{l-1}(A)$ is defined by
$$\lambda(f)=diag(\underbrace{f(x_{0}),f(x_{0}),\cdots,f(x_{0})}\limits_{l-1}),$$
for a fixed base point $x_{0}\in X$.\\

Let $A=\bigoplus_{i=1}^{{n}}M_{k_{i}}(C(X_{i}))$, where $X_{i}$ are connected finite simplicial complexes. A unital $\ast$-monomorphism $\phi:A\longrightarrow M_{l}(A)$ is called a (unital) simple embedding, if $\phi$ is of the form $\phi=\oplus\phi^{i}$ defined by
$$\phi(f_{1},f_{2},\cdots,f_{n})=(\phi^{1}(f_{1}),\phi^{2}(f_{2}),\cdots,\phi^{n}(f_{n})),$$
where the homomorphisms $\phi^{i}: A^{i}(=M_{k_{i}}(C(X_{i})))\longrightarrow M_{l}(A^{i})$ are unital simple embeddings.


\vspace{3mm}

\noindent\textbf{\S3. Additional decomposition theorems and factorization theorems}

\vspace{-2.4mm}

In this section, we will prove certain decomposition theorems which say that  the map $\phi_{n,m}$ (in an $AH$ inductive limit with the ideal property) can be decomposed  into two parts roughly described as below:\\
(a) the major part factors through as
$$A_{n}\xrightarrow{\xi_{1}}C\xrightarrow{\xi_{2}}A_{m},$$
where $C$ is a direct sum of matrix algebras over interval [0,1] or $\{pt\}$ and $\xi_{1},\xi_{2}$ are $\ast$ homomorphisms, and\\
(b) the other part factors through as
$$A_{n}\xrightarrow{\beta}B\xrightarrow{\alpha}A_{m},$$
where $B$ is a direct sum of matrix algebras over $C(S^{1}), C( T_{\uppercase\expandafter{\romannumeral2}, k}), C( T_{\uppercase\expandafter{\romannumeral3}, k})$ or $C(S^{2})$, $\alpha$ is a $\ast$ homomorphisms, but $\beta$ is a sufficiently multiplicative completely positive contraction. These theorems will play the roles of Theorem $5.32a$ and  Theorem $5.32b$ in [G4] in our proof. Since our limit algebra $A$ is no longer simple, we can not make each partial map to satisfy the dichotomy condition: each partial map $\phi^{i,j}_{n,m}$ is either injective or has finite dimensional image.

\noindent\textbf{3.1.} Recall that the total K-theory of $A$ is defined by
$$\underline{K}(A)=K_{\ast}(A)\oplus\bigoplus_{n=2}^{\infty}K_{\ast}(A,\mathbb{Z}/n).$$
Any $KK$ element $\alpha\in KK(A,B)$ determine a homomorphism
$$\underline{\alpha}: \underline{K}(A)\longrightarrow\underline{K}(B).$$
Now we let $W_{k}=T_{II,k}$ and let ${\mathcal{P}}\subset\mathop{\cup}\limits_{k}M_{\bullet}(A\otimes C(W_{k}\times S^{1}))$ be any finite set of projections. Then each element $p\in{\mathcal{P}}$ defines an element $[p]\in\underline{K}(A)$. We will use ${\mathcal{P}}\underline{K}(A)$ to denote the subgroup of $\underline{K}(A)$ generated by $\{[p]\in\underline{K}(A), p\in{\mathcal{P}}\}$. For each finite set $\mathcal{P}$, there is a finite set $G(\mathcal{P})\subset A$ (large enough) and $\delta(\mathcal{P})>0$ (small enough) such that if $\phi: A\longrightarrow B$ is $G(\mathcal{P})-\delta(\mathcal{P})$ multiplicative completely  positive contraction, then $\phi$ induces
$$\phi_{\ast}: \mathcal{P}\underline{K}(A)\longrightarrow \underline{K}(B).$$
(See [GL].)

Let $A=\bigoplus_{i=1}^{t_{n}}M_{[n,i]}(C(X_{n,i}))$, where $X_{n,i}$ are connected finite simplicial complexes. Then $K_{\ast}(A)$ is finitely generated, and there is a finite set $\mathcal{P}\subset\mathop{\cup}\limits_{k}M_{\bullet}(A\otimes C(W_{k}\times S^{1}))$ such that if two element $\alpha, \beta\in KK(A,B)$ satisfy
$$\underline{\alpha}\mid_{\mathcal{P}\underline{K}(A)}=\underline{\beta}\mid_{\mathcal{P}\underline{K}(A)}$$
then $\alpha=\beta\in KK(A,B)$.

The following is 5.24 of [G4].

\noindent\textbf{Definition 3.2.} For any finite set of projections $\mathcal{P}\subset\mathop{\cup}\limits_{k=2}\limits^{\infty}M_{\bullet}(A\otimes C(W_{k}\times S^{1}))$, let $G(\mathcal{P}), \delta(\cal{P})$ be as in 3.1. A $G(\mathcal{P})-\delta(\cal{P})$ multiplicative map $\phi:A\longrightarrow B$ is called quasi-$\mathcal{P}\underline{K}$-homomorphism  if there is a homomorphism $\psi:A\longrightarrow B$ with $\phi(\textbf{1}_{A})=\psi(\textbf{1}_{A})$ such that $[\phi]_{\ast}=[\psi]_{\ast}: \mathcal{P}\underline{K}(A)\longrightarrow \underline{K}(B)$.

\noindent\textbf{3.3.} Let $B=M_{\bullet}(C(Y))$, $Y$ is one of the space $T_{\uppercase\expandafter{\romannumeral2}, k}, T_{\uppercase\expandafter{\romannumeral3}, k}, S^{2}$, let $\mathcal{P}\subset\mathop{\cup}\limits_{k=2}\limits^{\infty} M_{\bullet}(B\otimes C(W_{k}\times S^{1}))$ be a finite set as in the end of 3.1 (one can choose $\mathcal{P}$ as in 5.16 of [G4]). Let $\phi:B\longrightarrow A$ be a homomorphism and let $\mathcal{P}_{1}\subset\mathop{\cup}\limits_{k=2}\limits^{\infty} M_{\bullet}(A\otimes C(W_{k}\times S^{1}))$, where $A$ and $\mathcal{P}_{1}$ satisfy the condition in 3.1. Furthermore we assume $\mathcal{P}_{1}\supseteq\phi(\mathcal{P})$. Then there is a finite set $G\subset A$ and $\delta>0$, such that if a $G-\delta$ multiplicative map $\psi: A\longrightarrow A_{1}$ is a quasi-$\mathcal{P}_{1}\underline{K}$-homomorphism, then $\psi\circ\phi:B\rightarrow A_{1}$ is a quasi-$\mathcal{P}\underline{K}$-homomorphism.

\noindent\textbf{Lemma 3.4.} Let $A=\bigoplus_{i=1}^{t_{n}}M_{[n,i]}(C(X_{n,i}))$. For any $\mathcal{P}_{1}\subseteq\mathop{\cup}\limits^{+\infty}_{k=2}M_{\bullet}(A\otimes C(W_{k}\times S^{1})),$ there is a finite set $F\subset A$ and $\varepsilon>0$, such that the following is true.

If $\phi: A\longrightarrow A^{\prime}=PM_{\bullet}(C(X))P$ is a unital homomorphism, $Q_{1},Q_{2}\in PM_{\bullet}(C(X))P$ are two projections with $Q_{1}+Q_{2}=P$, and $\phi_{1}\in Map(A, Q_{1}A^{\prime}Q_{1})_{1}$, $\phi_{2}\in Hom(A, Q_{2}A^{\prime}Q_{2})_{1}$ satisfy\\
(i) $\parallel \phi(f)-(\phi_{1}\oplus\phi_{2}(f))\parallel<\varepsilon$~~~~$\forall f\in F$ and\\
(ii) for each $i$, $\phi^{i}_{1}\in Map(A^{i}, Q_{1}A^{\prime}Q_{1})$ is either zero, or $$rank(\phi^{i}_{1}(\textbf{1}_{A^{i}}))\geqslant3(dim(X)+1)rank(\textbf{1}_{A^{i}})$$\\
then $\phi_{1}$ is a quasi-$\mathcal{P}_{1}\underline{K}$-homomorphism, where $A^{i}=M_{[n,i]}(C(X_{n,i}))$.
\begin{proof}
First, by Lemma 4.40 of [G4], if $F$ is large enough and $\varepsilon>0$ is small enough, then the condition (i) implies that $\phi_{1}$ is $G(\mathcal{P}_{1})-\delta(P_{1})$ multiplicative and induces maps on $\mathcal{P}_{1}\underline{K}(A)$.

We can assume that $A$  has only one block $M_{[n,1]}(C(X_{n,1}))$. Using Lemma 1.6.8 of [G4], one can assume that $\phi_1|_{M_{[n,1]}(\mathbb{C})}$ is a homomorphism, where $M_{[n,1]}(\mathbb{C})\subset M_{[n,1]}(C(X_{n,1}))$. Then one can reduce the proof of the lemma to the case that $[n,1]=1$, that is, $A=C(X_{n,1})$.
Then both $\phi\in Hom(C(X_{n,1}), PM_{\bullet}C(X)P)$ and
$\phi_{2}\in Hom(C(X_{n,1}), Q_{2}M_{\bullet}(C(X))Q_{2})$ induce $[\phi]\in kk(X,X_{n,1})$ and
$[\phi_{2}]\in kk(X,X_{n,1})$ (see 2.8 of [EG2], also see [DN] for notation $kk$). Since $kk(X,X_{n,1})$ is an Abelian group, and $rank (Q_{1})\geqslant3(dim(X)+1)rank(\textbf{1}_{A^{i}})$, by 6.4.4 of [DN] (see also Proposition 3.16 of [EG2]), there is a unital homomorphism $\psi_{1}\in Hom(A, Q_{1}A^{\prime}Q_{1})_{1}$ such that $$[\psi_{1}]=[\phi]-[\phi_{2}]\in kk(X,X_{n,1}).$$ Hence $[\psi_{1}]\mid_{\mathcal{P}_{1}\underline{K}(A)}=[\phi_{1}]\mid_{\mathcal{P}_{1}\underline{K}(A)}$, which implies $\phi_{1}$ is a quasi-$\mathcal{P}_{1}\underline{K}$-homomorphism.\\
\end{proof}
\noindent\textbf{Lemma 3.5.}~~Let $A=PM_{\bullet}(C(X))P$, where $X$ is connected finite simplicial complex. Let $F\subset A$ be approximately constant to within $\varepsilon$ (i.e. $\omega(F)<\varepsilon$). Then for any two homomorphisms $\phi,\psi: A\longrightarrow B=QM_{l}C(Y)Q$ defined by point evaluations with $K_{0}\phi=K_{0}\psi$ and  assuming that  for any $p\in A$,  $rank(\phi(p))\geqslant rank(p)\cdot dim(Y)$,  there exists a unitary $u\in B$ such that
$$\parallel\phi(f)-u\psi(f)u^{\ast}\parallel<2\varepsilon,~~~~~\forall f\in F.$$
\begin{proof}
Let $x_{0}\in X$ be a base point of $X$. There are finitely many points $\{x_{1},x_{2},\cdots,x_{n}\}\subset X$ such that $\phi$ factors through as
$$A\xrightarrow{\pi}A\mid_{\{x_{1},x_{2},\cdots,x_{n}\}}=\bigoplus_{i=1}^{n}M_{rank(P)}(\mathbb{C})\xrightarrow{\phi^{''}}B.$$
But for each $i$, there is a unitary $u_{i}\in M_{rank(P)}(\mathbb{C})$ such that
$$\parallel f(x_{0})-u_{i}f(x_{i})u^{\ast}_{i}\parallel<\varepsilon~~~~\forall f\in F,$$
since $\omega(F)<\varepsilon$. Let $\pi^{\prime}: A\longrightarrow \bigoplus_{i=1}^{n}M_{rank(P)}(\mathbb{C})$ defined by $$\pi^{\prime}(f)=(u^{\ast}_{1}f(x_{0})u_{1},u^{\ast}_{2}f(x_{0})u_{2},\cdots,u^{\ast}_{n}f(x_{0})u_{n}).$$ Then $\parallel\pi(f)-\pi^{\prime}(f)\parallel<\varepsilon$ for all $f\in F$. Evidently there is a homomorphism $\phi^{\prime}: M_{rank(P)}(\mathbb{C})\longrightarrow B$ such that $\parallel\phi(f)-\phi^{\prime}(f(x_{0}))\parallel<\varepsilon ,~\forall f\in F$. Similarly there is a $\psi^{\prime}: M_{rank(P)}(\mathbb{C})\longrightarrow B$ with $\parallel\psi(f)-\psi^{\prime}(f(x_{0}))\parallel<\varepsilon,~\forall f\in F$. On the other hand \\
$K_{0}\phi=K_{0}\psi$,  since $e_{\ast}: K_{0}(P M_{\bullet}(C(X))P)\longrightarrow K_{0}(M_{rank}(\mathbb{C}))$ is surjective for \\
$e: P M_{\bullet}(C(X))P\longrightarrow P(x_{0})M_{\bullet}(\mathbb{C})P(x_{0})\cong M_{rank(P)}(\mathbb{C})$. Furthermore by the condition that $rank(\phi(P))\geq rank(P)\cdot dim(Y)$ and the fact that for any two projections $Q_{1},Q_{2}$ (of rank at least $dim(Y)$) over $Y$, $[Q_{1}]=[Q_{2}]\in K_{0}(M_{\bullet}(C(Y)))$ implies $Q_{1}$ unitarily equivalent to $Q_{2}$, we know that there is a unitary $u\in B$ such that $\phi^{\prime}=u\psi^{\prime}u^{\ast}$.\\
\end{proof}
\noindent\textbf{Lemma 3.6.} Fix a positive integer $M\geq3$, suppose that $B=\bigoplus_{i=1}^{s}M_{l_{i}}(C(Y_{i}))$, where $Y_{i}$ are the spaces: $\{pt\}, [0,1], S^{1}, T_{\uppercase\expandafter{\romannumeral2}, k}, T_{\uppercase\expandafter{\romannumeral3}, k}, S^{2}$. Let $\varepsilon>0$ and  $$ \widetilde{G}(=\bigoplus\widetilde{G}^{i})\subset G(=\bigoplus G^{i})\subset \bigoplus B^{i}$$
be two finite sets satisfying that if $Y_{i}$ is one of $\{T_{\uppercase\expandafter{\romannumeral2}, k}\}_{k=1}^{\infty}, \{T_{\uppercase\expandafter{\romannumeral3}, k}\}_{k=1}^{\infty}$ and $S^{2}$, then $\omega(\widetilde{G}^{i})<\varepsilon$, and if $Y_{i}$ is one of $\{pt\}$, [0,1], $S^{1}$, then $\widetilde{G}^{i}=G^{i}$. Then there is a subset $G_{1}\subset B$ with $G_{1}\supset G(\mathcal{P})$, ($\mathcal{P}\subset\mathop{\cup}\limits^{\infty}\limits_{k=2}M_{\bullet}(B\otimes C(W_{k}\times S^{1}))$ as in the end of 3.1) and $\delta_{1}>0$ with $\delta_{1}<\delta(\mathcal{P})$, and a positive integer $L>0$ such that the following is true.

If a map $\alpha=\alpha_{0}\oplus\alpha_{1}: B\longrightarrow A=\bigoplus_{j=1}^{t}M_{k_{j}}(C(X_{j}))$, where $X_{j}$ are connected finite simplicial  complexes with $dim(X_{j})\leqslant M$, satisfying the following conditions:\\
(1) $\alpha_{0}$ is $G_{1}-\delta_{1}$ multiplicative, $\{\alpha_{0}(\textbf{1}_{B^{i}})\}^{s}_{i=1}$ are mutually othogonal projections, and for any block $B^{i}$ with $Y_{i}=T_{\uppercase\expandafter{\romannumeral2}, k}, T_{\uppercase\expandafter{\romannumeral3}, k}, S^{2}$ and any block $A^{j}$, the partial map $\alpha^{i,j}_{0}$ is quasi-$\mathcal{P}\underline{K}$-homomorphism; and\\
(2) $\alpha_{1}$ is a homomorphism defined by point evaluations and for each block $B^{i}$, with $Y_{i}=T_{\uppercase\expandafter{\romannumeral2}, k}, T_{\uppercase\expandafter{\romannumeral3}, k}, S^{2}$ and any block $A^{j}$, $$\alpha^{i,j}_{1}(\textbf{1}_{B_{i}})\geqslant L\alpha^{i,j}_{0}(\textbf{1}_{B_{i}}),$$ then there is a unital homomorphism $\alpha^{\prime}: B\longrightarrow \alpha(\textbf{1}_{B}) A\alpha (\textbf{1}_{B})$ such that $$\parallel\alpha^{\prime}(g)-\alpha(g)\parallel<3\varepsilon~~~\forall g\in\widetilde{G}.$$
\begin{proof}
Without loss of generality, we can assume $B$ has only one block: $B=M_{l}(C(Y))$. If $Y=\{pt\},[0,1],S^{1}$, the Lemma is true since $B$ is stably generated (see Lemma 1.6.1 of [G4]).

Suppose that $Y=T_{\uppercase\expandafter{\romannumeral2}, k}, T_{\uppercase\expandafter{\romannumeral3}, k}$ or $S^{2}$. Let $G_{1}$ and $\delta_{1}$, and $\eta$ be as Lemma 5.30 of [G4] (note that, for a positive integer $L$ and positive number $\eta$, the notion $PE(L,\eta)$ (used in Lemma 5.30 of [G4]) is defined in Definition 4.38 of [G4]). Let $\{y_{j}\}^{K}_{j=1}$ be an $\eta$-dense subset of $Y$. Let $L=4KM$. Let us verify the conclusion holds for such a choice. If $\alpha_{0}=0$, then one can choose $\alpha^{\prime}=\alpha_{1}$ as desired. Suppose $\alpha_{0}\neq0$ and therefore $rank$  $\alpha_{0}(\textbf{1}_{B})>0$. Using $ rank$ $\alpha_{1}(\textbf{1}_{B})\geqslant L rank(\alpha_{0}
(\textbf{1}_{B}))$, there is a unital
homomorphism $\alpha^{\prime}_{1}: B\longrightarrow \alpha_{1}(\textbf{1}_{B})A\alpha_{1}(\textbf{1}_{B})$ which satisfies the following conditions\\
(i) $\alpha^{\prime}_{1}$ is homotopy to $\alpha_{1}$\\
(ii) $\alpha^{\prime}_{1}$ is defined by direct sum of point evaluations at different points\\
(iii) all the points in the $\eta$-dense set $\{y_{j}\}^{K}_{j=1}$ are among those point evaluations that define $\alpha^{\prime}_{1}$ and each point evaluation at $y_{i}$ satisfies that  the  rank of the  image of $\textbf{1}_{B}$ is  at least rank $\alpha_{0}(\textbf{1}_{B})$---that is $\alpha^{\prime}_{1}$ has property $PE(rank$ $\alpha_{0}(\textbf{1}_{B}), \eta)$.

By 5.30 of [G4], there is a homomorphism $\alpha{''}: B\longrightarrow \alpha(\textbf{1}_{B})A\alpha(\textbf{1}_{B})$ such that $$\parallel\alpha{''}(f)-(\alpha_{0}\oplus\alpha^{\prime}_{1}(f))\parallel<\varepsilon,~~~\forall f\in \widetilde{G}.$$ On the other hand by Lemma 3.5, there is a unitary $u\in \alpha^{\prime}_{1}(\textbf{1}_{B})
A\alpha^{\prime}_{1}(\textbf{1}_{B})$ (note that $\alpha^{\prime}_{1}(\textbf{1}_{B})=\alpha_{1}(\textbf{1}_{B})$) such that
$$\parallel u^{\ast}\alpha^{\prime}_{1}(f)u-\alpha_{1}(f)\parallel<2\varepsilon, ~~~\forall f\in\widetilde{G}$$
Let $\alpha^{\prime}(t)=diag(\alpha_{0}(\textbf{1}_{B}),u^{\ast})\alpha{''}(f)diag(\alpha_{0}(\textbf{1}_{B}),u)$.\\
Evidently, we have
 $$\parallel\alpha{'}(f)-(\alpha_{0}\oplus\alpha_{1})(f)\parallel<3\varepsilon, ~~~\forall f\in \widetilde{G}.$$
 \end{proof}
\noindent\textbf{Lemma 3.7.} Let $M$ be a fixed positive integer. Let $B=M_{l}(C(Y))$, $Y=T_{\uppercase\expandafter{\romannumeral2}, k}, T_{\uppercase\expandafter{\romannumeral3}, k}$ or $S^{2}$. Let the set of projection $\mathcal{P}\subset M_{\bullet}(B\otimes C(W_{k}\times S^{1}))$ be as in 3.1 (also see 5.16 of [G4]).

Let $A=RM_{l}(C(X))R$ with $X$ a connected finite simplicial complex and let $\alpha: B\longrightarrow A$ be a homomorphism. Let $\mathcal{P}^{\prime}\subset\mathop{\cup}\limits^{\infty}\limits_{k=2}M_{\bullet}(A\otimes C(W_{k}\times S^{1}))$ be a set of projections (chosen for  $A$)  as in the end of 3.1 and $\mathcal{P}^{\prime}\supseteq(\alpha\otimes id)(\mathcal{P})$.

For any finite sets $\widetilde{G}_{1}\subset G_{1}\subset B$ and numbers $\varepsilon_{1}>0, \delta_{1}>0,$ with $\omega(\widetilde{G}_{1})<\varepsilon_{1},$ there are a finite set $G_{2}\subset A$ a number $\delta_{2}>0$, and positive integer $L$, such that the following are true. Let $C=M_{\bullet}(C(Z))$ with $Z$ connected simplicial complex and $dim(Z)\leqslant M$, and let $Q_{0},Q_{1}\in C$ be two orthogonal projections.\\
(1) If $\psi_{0}: A\longrightarrow Q_{0}CQ_{0}$ is $G_{2}-\delta_{2}$ multiplicative quasi- $\mathcal{P}^{\prime}\underline{K}$-homomorphism and $\psi_{0}(\alpha(\textbf{1}_{B}))$ is a projection, then $\psi_{0}\circ\alpha$ is a $G_{1}-\delta_{1}$ multiplicative quasi-$\mathcal{P}\underline{K}$-homomorphism.\\
(2) If $\psi_{0}$ as in (1) and $\psi_{1}: A\longrightarrow Q_{1}CQ_{1}$ is defined by point evaluations with\\
 rank $(\psi_{1}(\textbf{1}_{A}))\geq L\mbox{rank}(\psi_{0}(\textbf{1}_{A}))$, then there is a homomorphism \\
 $\psi: B\longrightarrow(Q_{0}\oplus Q_{1})C(Q_{0}\oplus Q_{1})$ such that
$$\parallel \psi(g)-(\psi_{0}\oplus \psi_{1})(\alpha(g))\parallel<3\varepsilon.$$
\begin{proof}
The proof is the same as the proof of Lemma 5.31 of [G4] using Lemma 3.6 above to replace 5.30 of [G4].\\
\end{proof}

\noindent\textbf{Theorem 3.8.} Let $M\geqslant3$ be a positive integer. Let $\lim\limits_{n\rightarrow\infty}\limits^{}(A_{n}=\bigoplus_{i=1}^{k_{n}}M_{[n,i]}(C(X_{n,i})), \phi_{n,m})$ be an $AH$ inductive limit such that the limit algebra has the ideal property, where $X_{n,i}$ are connected simplicial complexes with $dim(X_{n,i})\leq M$, for all $n,i$. Let $B=\bigoplus_{i=1}^{s}M_{l_{i}}(C(Y_{i})),$ where $Y_{i}$ are spaces $\{pt\}$, [0,1], $S^{1}$, $T_{\uppercase\expandafter{\romannumeral2}, k}$, $T_{\uppercase\expandafter{\romannumeral3}, k},$ and $S^{2}$. Suppose that $\widetilde{G}(=\bigoplus\widetilde{G}^{i})\subset G(=\bigoplus G^{i})\subset B(=\bigoplus B^{i})$ are two finite sets and $\varepsilon_{1}>0$ is a positive number with $\omega(\widetilde{G}^{i})<\varepsilon_{1}$, if $Y_{i}=T_{\uppercase\expandafter{\romannumeral2}, k}, T_{\uppercase\expandafter{\romannumeral3}, k},$ or $S^{2}$. Suppose that $L$ is any positive integer. Let $\alpha: B\longrightarrow A_{n}$ be any homomorphism. Denote $\alpha(\textbf{1}_{B}):= R(=\bigoplus R^{i})\in A_{n}(=\bigoplus A^{i}_{n})$. Let $F\subset RA_{n}R$ be any finite set and $\varepsilon<\varepsilon_{1}$ be any positive number.

It follows that there are $A_{m}$, and mutually orthogonal projections $Q_{0},Q_{1},Q_{2}\in A_{m}$ with $\phi_{n,m}(R)=Q_{0}+Q_{1}+Q_{2},$ a unital map $\theta_{0}\in Map(RA_{n}R, Q_{0}A_{m}Q_{0})_{1}$, two unital homomorphisms $\theta_{1}\in Hom(RA_{n}R, Q_{1}A_{m}Q_{1})_{1}$, $\xi\in Hom(RA_{n}R, Q_{2}A_{m}Q_{2})_{1}$ such that\\
(1) $\parallel \phi_{n,m}(f)-(\theta_{0}(f)\oplus\theta_{1}(f)\oplus \xi(f))\parallel<\varepsilon~~~~\forall f\in F$\\
(2) there is a homomorphism $\alpha_{1}: B\longrightarrow(Q_{0}+Q_{1})A_{m}(Q_{0}+Q_{1})$ such that
\begin{center}
$\parallel \alpha_{1}(g)-(\theta_{0}\oplus \theta_{1})\circ\alpha(g)\parallel<3\varepsilon_{1}~~~~\forall g\in \widetilde{G}^{i},
~~if\; Y_{i}= T_{\uppercase\expandafter{\romannumeral2}, k}, T_{\uppercase\expandafter{\romannumeral3}, k}$ or $ S^{2}$
\end{center}
and
\begin{center}
$\parallel \alpha_{1}(g)-(\theta_{0}\oplus\theta_{1})\circ\alpha(g)\parallel<\varepsilon~~~~\forall g\in G^{i},~~ if\; Y_{i}=\{pt\}, [0,1]$ or $S^{1}.$\\
\end{center}
(3) $\theta_{0}$ is $F-\varepsilon$ multiplicative and $\theta_{1}$ satisfies that for any nonzero projections $e\in R^{i}A^{i}_{n}R^{i}$
$$  \theta^{i,j}_{1}([e])\geqslant L\cdot[\theta^{i,j}_{0}(R^{i})].$$
(4) $\xi$ factors through a $C^{\ast}$-algebra $C$---a direct sum of matrix algebras over C[0,1] or $\mathbb{C}$ as
$$\xi:RA_{n}R\xrightarrow{\xi_{1}}C\xrightarrow{\xi_{2}}Q_{2}A_{m}Q_{2}.$$
\begin{proof}
Let $D\subset A_{n}=\oplus A^{i}_{n}$ be defined by
$$D=\bigoplus\limits_{j}(\bigoplus\limits_{i}\alpha^{i,j}(\mathbb{C}\cdot\textbf{1}_{B^{i}}))\subset \bigoplus\limits_{j}A^{j}_{n}$$
which is a finite dimensional subalgebra of $A_{n}$ containing the set of mutually orthogonal projections $\{E^{i,j}=\alpha^{i,j}(\textbf{1}_{B^{i}})\}_{i,j}$.

Apply Corollary 2.22 for sufficiently large set $F^{\prime}\subset RA_{n}R$, sufficiently small number $\varepsilon^{\prime}>0$ and sufficiently large integer $J>0$, to obtain $A_{m}$ and the decomposition $\theta_{0}\oplus\theta_{1}\oplus\xi$ of $\phi_{n,m}\mid_{RA_{n}R}$ as $\psi^{\prime}_{0}\oplus\psi^{\prime}_{1}\oplus\psi^{\prime}_{2}$ in the corollary. By Lemma 1.6.8 of [G4], we can assume $\theta_{0}\mid_{D}$ is a homomorphism. The condition (1) follows if we choose $F^{\prime}\supset F$, and $\varepsilon^{\prime}<\varepsilon$. The $F-\varepsilon$ multiplicative of $\theta_{0}$ in (3) follows from Lemma 4.40 of [G4], if $F^{\prime}$ is large enough and $\varepsilon^{\prime}$ is small enough; and property of $\theta_{1}$ in (3) follows if we choose $J>L$.

To construct $\alpha_{1}$ as desired in the condition (2), we need to construct
$$\alpha^{i,j,k}_{1}: B^{i}\longrightarrow\theta^{j,k}(E^{i,j})A^{k}_{m}\theta^{j,k}(E^{i,j})$$
where $\theta=\theta_{0}\oplus \theta_{1}$, to satisfy
$$\parallel \alpha^{i,j,k}_{1}(g)-\theta^{j,k}\circ\alpha^{i,j}(g)\parallel<3\varepsilon_{1}~~~~\forall g\in \widetilde{G}^{i},
~~if\; Y_{i}= T_{\uppercase\expandafter{\romannumeral2}, k}, T_{\uppercase\expandafter{\romannumeral3}, k}, S^{2}$$
and
$$\parallel \alpha^{i,j,k}_{1}(g)-\theta^{j,k}\circ\alpha^{i,j}(g)\parallel<\varepsilon~~~~\forall g\in {G}^{i},~~ if\; Y_{i}=\{pt\}, [0,1], S^{1}.$$\\
For the case of $Y_{i}=\{pt\}, [0,1], S^{1}$, the existence of $\alpha^{i,j,k}_{1}$ follows from Lemma 1.6.1 and Lemma 4.40 of [G4]---that is $\theta^{j,k}_{0}\circ\alpha^{i,j}$ itself can be perturbed to a homomorphism as $B^{i}$ is stably generated.\\
For the case $ Y_{i}= T_{\uppercase\expandafter{\romannumeral2}, k}, T_{\uppercase\expandafter{\romannumeral3}, k}, S^{2}$, the existence of $\alpha^{i,j,k}_{1}$ follows from Lemma 3.7 and $\omega(\widetilde{G}^{i})<\varepsilon_{1}$ provided $J$ is large enough, $F^{\prime}$ large enough, $\varepsilon^{\prime}$ small enough.\\
\end{proof}
\noindent\textbf{Theorem 3.9.}  Let $M$ be a positive integer. Let $\lim\limits_{n\rightarrow\infty}\limits^{}(A_{n}=\bigoplus_{i=1}^{k_{n}}M_{[n,i]}(C(X_{n,i})), \phi_{n,m})$ be an $AH$ inductive limit such that the limit algebra has the ideal property, where $X_{n,i}$ are connected simplicial complexes with $dim(X_{n,i})\leq M$, for all $n,i$. Let $B=\bigoplus_{i=1}^{s}M_{l_{i}}(C(Y_{i})),$ where $Y_{i}$ are spaces $\{pt\}$, [0,1], $S^{1}$, $T_{\uppercase\expandafter{\romannumeral2}, k}$, $T_{\uppercase\expandafter{\romannumeral3}, k}$ and $S^{2}$. Suppose that $\widetilde{G}(=\bigoplus\widetilde{G}^{i})\subset G(=\bigoplus G^{i}\subset B(=\bigoplus B^{i})$ are two finite subsets and $\varepsilon_{1}$ is a positive number such that $\omega(\widetilde{G}^{i})<\varepsilon_{1}$, if $Y_{i}$ is one of $T_{\uppercase\expandafter{\romannumeral2}, k}, T_{\uppercase\expandafter{\romannumeral3}, k},$ or  $ S^{2}$. Let $\alpha: B\longrightarrow A_{n}$ be a homomorphism and $F(\supset\alpha(G))$ be a finite subset of $A_{n}$ and $\varepsilon<\varepsilon_{1}$ be any positive number.

It follows that there are $A_{m}$ and mutually orthogonal projections $P,Q\in A_{m}$ with $\phi_{n,m}(\textbf{1}_{A_{n}})=P+Q$, a unital map $\theta\in Map(A_{n}, PA_{m}P)_{1}$, and a unital homomorphism $\xi\in Hom(A_{n}, QA_{m}Q)_{1}$ such that\\
(1) $\parallel \phi_{n,m}(f)-(\theta(f)\oplus\xi(f))\parallel<\varepsilon~~~~\forall f\in F$\\
(2) there is a homomorphism $\alpha_{1}: B\longrightarrow PA_{m}P$ such that
$$\parallel \alpha^{i,j}_{1}(g)-(\theta\circ\alpha)^{i,j}(g)\parallel<3\varepsilon_{1}~~~~\forall g\in \widetilde{G}^{i},
~~if\; Y_{i}= T_{\uppercase\expandafter{\romannumeral2}, k}, T_{\uppercase\expandafter{\romannumeral3}, k}, S^{2},$$ $$\parallel \alpha^{i,j}_{1}(g)-(\theta\circ\alpha)^{i,j}(g)\parallel<\varepsilon~~~~\forall g\in G^{i},~~ if\; Y_{i}=\{pt\}, [0,1], S^{1}.$$\\
(3) $\omega(\theta(F))<\varepsilon$ and $\theta$ is $F-\varepsilon$ multiplicative.\\
(4) $\xi$ factors through a $C^{\ast}$-algebra $C$---a direct sum of matrix algebras over $C[0,1]$ or $\mathbb{C}$ as
$$\xi: A_{n}\xrightarrow{\xi_{1}}C\xrightarrow{\xi_{2}}QA_{m}Q.$$\\
The proof is similar to but slightly easier than that of Theorem 3.8, and we omit it.

The following is Corollary 1.6.15 of [G4].

\noindent\textbf{Proposition 3.10.} Let $A=\bigoplus_{k=1}^{l}M_{s(k)}(C(X_{k}))$, where $X_{k}$ are connected finite simplicial complexes and $s(k)$ are positive integers. Let $F\subset A$ be a finite set and $\varepsilon>0$. There are a finite set $G\subset A$ and a number $\delta>0$ with the following property. If $B$ is a unital $C^{\ast}$-algebra, $p\in B$ is a projection, $\phi_{t}: A\longrightarrow pBp, 0\leqslant t\leqslant1$ is a continuous path of $G-\delta$ multiplicative maps, then there are a positive integer $L$, and $\eta>0$ such that for a homomorphism $\lambda: A\longrightarrow B\otimes \mathcal{K}$ (with finite dimensional image), there is a unitary $u\in(p\oplus\lambda(\textbf{1}))(B\otimes \mathcal{K})(p\oplus\lambda(\textbf{1}))$ satisfying
$$\parallel \phi_{0}(f)\oplus\lambda(f)-u(\phi_{1}(f)\oplus\lambda(f))u^{\ast}\parallel<\varepsilon,~~~~\forall f\in F,$$
provided that $\lambda$ is of the following form: there are an $\eta$-dense subset $\{x_{1},x_{2},\cdots,x_{\bullet}\}\subset \amalg^{l}_{k=1}X_{k}(=Sp(A))$, and a set of mutually orthogonal projections $$\{p_{1},p_{2},\cdots,p_{\bullet}\}\subset \lambda(\bigoplus_{k}e^{k}_{11})(B\otimes \mathcal {K})\lambda(\bigoplus_{k}e^{k}_{11})$$ with $[p_{i}]\geqslant L\cdot[p]$, such that
$$\lambda(f)=\sum\limits^{n}\limits_{i=1}p_{i}\otimes f(x_{i})\oplus\lambda_{1}(f),~~~~\forall f\in A,$$
where $\lambda_{1}$ is also a homomorphism, under the identification $$\lambda(\textbf{1}_{A^{k}})B\lambda(\textbf{1}_{A^{k}})\cong \lambda(e^{k}_{11})B\lambda(e^{k}_{11})\otimes M_{s(k)}(\mathbb{C}).$$

\noindent\textbf{Corollary 3.11.}~~Let $A=\bigoplus_{k=1}^{l}M_{s(k)}(C(X_{k}))$, where $X_{k}$ are connected finite simplicial complexes and $s(k)$ are positive integers. Let $\varepsilon>0$ and $F\subset A$ be a finite set with\\
$\omega(F)<\varepsilon$. There are a finite set $G\subset A$ and a number $\delta>0$ with the following property. If $B=\bigoplus\limits_{j=1}^{}M_{r(j)}(C(Z_{j}))$, where $Z_{j}$ are connected finite simplicial complexes and $p\in B$ is a projection, $\phi_{t}: A\longrightarrow pBp, 0\leqslant t\leqslant1$, is a continuous path of $G-\delta$ multiplicative maps, then there is a positive integer $L$ such that for a homomorphism $\lambda: A\longrightarrow B\otimes \mathcal{K}$ with finite dimensional image, there is a unitary $u\in(p\oplus\lambda(\textbf{1}))(B\otimes \mathcal{K})(p\oplus\lambda(\textbf{1}))$ satisfying
$$\parallel \phi_{0}(f)\oplus\lambda(f)-u(\phi_{1}(f)\oplus\lambda(f))u^{\ast}\parallel<5\varepsilon~~~~\forall f\in F,$$ provided $\lambda^{k}=\lambda\mid_{A^{k}}: M_{s(k)}(C(X_{k}))\longrightarrow B\otimes \mathcal{K}$ has finite dimensioned image (or equivalently, is defined by point evaluation) with
$$[\lambda^{k}(\textbf{1}_{A^{k}})]\geqslant L\cdot[p]\in K_{0}(B).$$
\begin{proof}
Let $L^{\prime}$ and $\eta$ be $L$ and $\eta$ as in Proposition 3.10. Let $\{x_{1},x_{2},\cdots,x_{K}\}\subset\mathop{\amalg}\limits^{l}\limits_{k=1}X_{k}$ be a $\eta$-dense subset. Choose $L=4L^{\prime}\cdot K \max(\mathop{dim}\limits_{j}(Z_{j})+1)$. If $\lambda$ satisfies the condition in our corollary then it is easy to find $\lambda^{\prime}: A\longrightarrow \lambda(\textbf{1})(B\otimes\mathcal{K})\lambda(\textbf{1})$ satisfies the condition in Proposition 3.10 and $K_{0}\lambda=K_{0}\lambda^{\prime}$. Then the corollary follows from Proposition 3.10 and Lemma 3.5.\\
\end{proof}
\noindent\textbf{Theorem 3.12.} Let $B_{1}=\bigoplus_{j=1}^{s}M_{k(j)}(C(Y_{j}))$, where $Y_{j}$ are spaces
$\{pt\}$, [0,1], $S^{1}$, $\{T_{\uppercase\expandafter{\romannumeral2}, k}\}^{\infty}_{k=2}$, $\{T_{\uppercase\expandafter{\romannumeral3}, k}\}^{\infty}_{k=2}$ and $S^{2}$. Let $X$ be a connected finite simplicial complex and let $A=M_{N}(C(X))$. Let $\widetilde{G}_{1}(=\bigoplus\widetilde{G}^{i}_{1})\subset G_{1}(=\bigoplus G^{i}_{1})\subset B_1(=\bigoplus B_1^{i})$ be two finite sets with $\omega(\widetilde{G}^{i}_{1})<\varepsilon_{1}$ for certain $\varepsilon_{1}>0$ and any $i$ with $Y_{i}$ being $T_{\uppercase\expandafter{\romannumeral2}, k}$, $T_{\uppercase\expandafter{\romannumeral3}, k}$ or $S^{2}$. Let $\alpha_{1}:B_{1}\rightarrow A$ be a homomorphism, and let $F_{1}\subset A$ be a finite set and any positive number $\varepsilon<\varepsilon_{1}$ and $\delta>0$. Then there exists a diagram
$$
\xymatrix{
 A\ar[rdrd]^{\beta}\ar[rr]^{\phi} &  &  A^{\prime}  \\
     &  &           &  &         &  &          \\
     B_{1}\ar[uu]_{\alpha_{1}}\ar[rr]^{\psi}  &  &  B_{2}\ar[uu]_{\alpha_{2}}     \\
 }
$$
where $A^{\prime}=M_{L}(A)$, and $B_{2}$ is a direct sum of matrix algebras over space: $\{pt\}$, [0,1], $S^{1}$, $T_{\uppercase\expandafter{\romannumeral2}, k}$, $T_{\uppercase\expandafter{\romannumeral3}, k}$ and $S^{2}$, with the following conditions \\
(1) $\psi$ is a homomorphism, $\alpha_{2}$ is a unital injective homomorphism and $\phi$ is a unital simple embedding (see 2.23 for the definition of unital simple embedding);\\
(2) $\beta\in Map(A,B_{2})_{1}$ is $F_{1}-\delta$ multiplicative;\\
(3) $\parallel\psi(g)-\beta\circ\alpha_{1}(g)\parallel<5\varepsilon_{1}~~~\forall g\in\widetilde{G}^{i}_{1}$,~~ if $Y_{i}= T_{\uppercase\expandafter{\romannumeral2}, k}, T_{\uppercase\expandafter{\romannumeral3}, k}, S^{2}$,\\
$~~~~~\parallel \psi(g)-\beta\circ\alpha_{1}(g)\parallel<\varepsilon~~~\forall g\in G^{i}_{1}$,~~ if  $Y_{i}=\{pt\}, [0,1], S^{1}$,\\
$~~~~~~ \mbox{and}$ \\$~~~~~\parallel \phi(f)-\alpha_{2}\circ\beta(f)\parallel<\varepsilon$~~~$\forall f\in F_{1}$; and\\
(4) $\beta(F_{1})\cup \psi(G_{1})\subset B_{2}$ satisfies that
    $~\omega(\beta(F_{1})\cup \psi(G_{1}))<\varepsilon.$
\begin{proof}
The proof is similar to the proof of Theorem 1.6.26 (see 1.6.25) of [G4]. The differences are the following. First, we do not have the condition that $\alpha_{1}$ is injective, so we need to use Corollary 3.11 to deal with blocks $B^{i}_{1}$ with $Y_{i}= T_{\uppercase\expandafter{\romannumeral2}, k}, T_{\uppercase\expandafter{\romannumeral3}, k}$ or $S^{2}$. For those block $B^{i}_{1}$ with $Y_{i}=\{pt\}, [0,1], S^{1}$, we use the condition that $B^{i}_{1}$ is stably generated (see [G4] 1.6.1). Secondly,  we need to make the condition (4) hold.

Without loss of generality, we can assume $\alpha_{1}(B^{i}_{1})\neq 0$ for each block  $B^{i}_{1}$, otherwise we simply take $\psi$ to be zero on this block.

Apply Corollary 3.11 to $\widetilde{G}^{i}_{1}\subset B^{i}_{1}$ for blocks $Y_{i}= T_{\uppercase\expandafter{\romannumeral2}, k}, T_{\uppercase\expandafter{\romannumeral3}, k}$ or $S^{2}$, to obtain $G^{i}\subset B^{i}_{1}$ and $\delta_{1}$ as in the corollary. Apply Lemma 1.6.1 of [G4] to $G^{i}_{1}\subset B^{i}_{1}$ for $Y_{i}=\{pt\}, [0,1], S^{1}$, to obtain $G^{i}$ and $\delta^{\prime}_{1}$. We assume $G^{i}\supset G_{1}^{i}$. Let $G=\bigoplus G^{i}$. Let $F^{\prime}_{1}\subset A$ and $\delta_{2}>0$ be such that if $\beta: A\longrightarrow C$ (any $C^{\ast}$-algebra) is $F^{\prime}_{1}-\delta_{2}$ multiplicative, then $\beta\circ\alpha_{1}$ is $G-\min(\delta_{1},\delta^{\prime}_{1})$ multiplicative.

Apply Proposition 3.10 to $F^{\prime}_{1}\cup F_{1}$ (as $F\subset A$) and $\min(\delta_{2},\delta)$ (as $\varepsilon>0$) to obtain $F\subset A$ and $\delta_{3}$ (in place of the set $G$ and the number $\delta$ there). We can assume $\delta_{3}<\min(\delta_{2},\delta,\delta_{1},\delta^{\prime}_{1})$ and $F\supset F^{\prime}_{1}\cup F_{1}$.

Apply 1.6.24 of [G4] to obtain the following diagram
$$
\xymatrix{
 A\ar[rdrd]^{\beta}\ar[rr]^{\phi} &  &  A^{\prime}  \\
     &  &           &  &         &  &          \\
     B_{1}\ar[uu]_{\alpha_{1}}\ar[rr]^{\psi}  &  & B_{2}\ar[uu]_{\alpha_{2}~~~~~~~~~~~~~~~~~~~~~~~~~~~~~~~~{\displaystyle\huge{(I)}}}
 }
$$
where $A^{\prime}=M_{L}(A)$, and $B_{2}$ is a direct sum of matrix algebras over the spaces $\{pt\}, [0,1], S^{1}$, $\{T_{\uppercase\expandafter{\romannumeral2}, k}\}, \{T_{\uppercase\expandafter{\romannumeral3}, k}\}$ and $S^{2}$ with the following conditions:\\
(i) $\psi$ is homomorphism, $\alpha_{2}$ is a unital injective homomorphism and $\phi$ is a unital simple embedding\\
(ii) $\beta\in Map(A,B_{2})_{1}$ is $F-\delta_{3}$ multiplicative and therefor $\beta\circ\alpha_{1}$ is $G-min(\delta_{1},\delta^{'}_{1})$ multiplicative\\
(iii) there exist two homotopies  $\Psi\in Map(B_{1},B_{2}[0,1])$ and $\Phi\in Map(A,A^{'}[0,1])$ such that $\Psi$ is $G-\delta_{3}$ multiplicative and $\Phi$ is $F-\delta_{3}$ multiplicative.\\
Note that we can choose a simple embedding $\xi:B_{2}\rightarrow M_{k}(B_{2})$ such that $\omega(\xi(\beta(F)\cup\psi(G)))<\varepsilon$ by adding to the identity the homomorphism defined by point evaluations at all points in a sufficiently dense finite set.
Let $\xi^{'}:A^{'}\rightarrow M_{k}(A^{'})$ be any simple embedding.
Since $\alpha_{2}$ take trivial projections to trivial projections(see Remark 1.6.20 of [G4]), we know that
$$
\xymatrix{
A^{'}\ar[rr]^{\xi^{'}}  &  &  M_{k}(A^{'}) &  &   \\
     &  &  &  &  &  &  \\
     B_{2}\ar[rr]^{\xi}\ar[uu]^{\alpha_{2}} &  &  M_{k}(B_{2})\ar[uu]_{\alpha_{2}\otimes id_{k}}
      &  &   \\
 }
$$
commutes up to homotopy. Therefor replacing  $\beta$, $\phi$, $\psi$, and $\alpha_2$, 
by $\xi\circ\beta$, $\xi^{'}\circ\phi$,  $\xi\circ\psi$, and $\alpha_{2}\otimes id_{k}$, respectively,  we can assume our original diagram
(I) also satisfies \\
(iv) $\omega(\beta(F)\cup\psi(G))<\varepsilon$.

Without loss of generality, we can assume $B_{1}$ has only one block $Y=T_{II,k}, T_{III,k}$ or $S^{2}$. Regarding the homotopy $\Psi$ as the homotopy path $\phi_{t}$ in Corollary 3.11, we can obtain $L_{1}$ as the number $L$ in Corollary 3.11.
Similarly regarding the homotopy $\Phi$ as $\phi_{t}$ in Proposition 3.10 we obtain $L_{2}$ (as $L$) and $\eta$.

Choose $\{x_{1},x_{2},\cdot\cdot\cdot,x_{m}\}\subset X$ to be a $\eta$-dense subset. Define $$\lambda_{1}:A(=M_{N}(C(X)))\rightarrow M_{mN}(B_{2})$$ by
$$\lambda_{1}(f)=diag(\textbf{1}_{B_{2}}\otimes f(x_{1}),\textbf{1}_{B_{2}}\otimes f(x_{2}),\cdot\cdot\cdot,\textbf{1}_{B_{2}}\otimes f(x_{m}))$$

Let $L{'}=max\{L_{1},L_{2}\}$ and $n=mNL{'}$. Define
$$ \lambda:A\longrightarrow M_{n}(B_{2})=M_{L{'}}(M_{mN}(B_{2}))\;\;by$$
$$ \lambda=diag(\underbrace{\lambda_{1},\lambda_{1},\cdot\cdot\cdot,\lambda_{1}}\limits_{L^{\prime}})$$\\
Then $\lambda\circ\alpha_{1}:B_{1}\longrightarrow M_{n}(B_{2})$ satisfies the condition for $\lambda$ in Corollary 3.11 for the homotopy $\Psi$ and the positive integer $L_1$. Also  $(\alpha_{2}\otimes id_{n})\circ \lambda:A\longrightarrow M_{n}(A^{\prime})$
satisfies the condition for $\lambda$ in Proposition 3.10 for the homotopy $\Phi, L_{2}$ and $\eta$. Therefor there are $u_{1}\in M_{n+1}(B_{2})$ and $u_{2}\in M_{n+1}(A^{\prime})$ such that
 $$\|(\beta\oplus\lambda)\circ\alpha_{1}(g)-u_{1}((\psi\oplus\lambda\circ\alpha_{1})(g))u^{*}_{1}\|\leq5\varepsilon_{1},~~~\forall g\in \widetilde{G_{1}} \eqno (\ast)$$\\
 and
 $$\|(\phi\oplus((\alpha_{2}\otimes id_{n})\circ \lambda))(f)-u_{2}((\alpha_{2}\otimes id_{n+1})\circ(\beta\oplus\lambda)(f))u^{*}_{2}\|<\varepsilon,~~~ \forall f\in F_{1}$$
 Note that if $Y_{i}=\{pt\}, [0,1]$ or $S^{1}$, then $\beta\circ\alpha_{1}$ itself is close to a homomorphism $\psi{'}:B_{1}\rightarrow B_{2}$. That is,  we can replace the above $(*)$ by
 $$\|\beta\circ\alpha_{1}(g)-\psi{'}(g)\|<\varepsilon,~~~ \forall g\in G_{1}.  \eqno (\ast\ast)   $$
 In the diagram (I) if we replace $B_{2}$ by $M_{n+1}(B_{2}), A^{\prime}$ by $M_{n+1}(A^{\prime}), \psi$ by $Adu_{1}\circ(\psi\oplus\lambda\circ\alpha_{1})$ (or $\psi$ by $\psi{'}\oplus(\lambda\circ\alpha_{1})$ for the case $Y_{i}=\{pt\}, [0,1],S^{1}$  using $(**)$), $\beta$ by $\beta\oplus\lambda, \alpha_{2}$ by $Adu_{2}\circ(\alpha_{2}\otimes id_{n+1})$ and finally $\phi$ by $\phi\oplus((\alpha_{2}\otimes id_{n})\circ\lambda)$, we get the desired diagram.\\
\end{proof}
 \noindent\textbf{3.13}~~Recall that in 1.1.7(h) of [G4], for $A=\bigoplus^{t}_{i=1}M_{k_{i}}(C(X_{i}))$,
 where $X_{i}$ are path connected simplicial complexes, we used the notation $r(A)$ to denote $\bigoplus^{t}_{i=1}M_{k_{i}}(\mathbb{C})$, which could be considered to be the subalgebra consisting of t-tuples
 of constant functions from $X_{i}$ to $M_{k_{i}}(\mathbb{C})(i=1,2,...,t)$.
 Fixed a base point $x^{0}_{i}\in X_{i}$ for each $X_{i}$, one defines a map $r:A\rightarrow r(A)$ by
 $$r(f_{1},f_{2},...,f_{t})=(f_{1}(x^{0}_{1}),f_{2}(x^{0}_{2}),\cdot\cdot\cdot,f_{t}(x^{0}_{t}))\in r(A)$$
 We have the following corollary.

 \noindent\textbf{Corollary 3.14}~~Let $B_{1}=\bigoplus^{s}_{j=1}M_{k(j)}(C(Y_{j}))$, where $Y_{j}$
 are the spaces $\{pt\}, [0,1], S^{1}$, $\{T_{II,k}\}_{k}$, $\{T_{III,k}\}_{k}$ and $S^{2}$.
 Let $A=\bigoplus^{t}_{j=1}M_{l(j)}(C(X_{j}))$, where $X_{j}$ are connected simplicial complexes. Let $\alpha_{1}:B_{1}\rightarrow A$ be any homomorphism. Let $\varepsilon_{1}>\varepsilon_{2}>0$ be any two positive numbers. Let $\widetilde{E}(=\bigoplus\widetilde{E}^{i})\subset E(=\bigoplus E^{i})\subset B_{1}(=\bigoplus B_1^{i})$ be two finite subsets with the condition
 $$\omega(\widetilde{E}^{i})<\varepsilon_{1}\;\; for\; all \;Y_{i}=T_{II,k}, T_{III,k}\; or\; S^{2}.$$
 Let $F\subset A$ be any finite subset, $\delta>0$. Then there exists a diagram
 $$
\xymatrix{
A\ar[rr]^{\phi\oplus r}\ar[rdrd]^{\beta\oplus r}  &  &  A^{'}\oplus r(A)  &  &    \\
     &  &  &  &  &  &  \\
     B_{1}\ar[rr]^{\psi\oplus(r\circ\alpha_{1})}\ar[uu]^{\alpha_{1}} &  &  B_{2}\oplus r(A)\ar[uu]_{\alpha_{2}\oplus id} &  &  \\
 }
$$
 where $A^{'}=M_{L}(A), B_{2}$ is a direct sum of matrix algebras over the spaces $\{pt\}, [0,1]$, $S^{1}$, $\{T_{II,k}\}, \{T_{III,k}\}$ and $S^{2}$, with the following properties.\\
(1) $\psi$ is a homomorphism, $\alpha_{2}$ is injective homomorphism and $\phi$ is a unital simple embedding (see 2.23).\\
(2) $\beta\in Map(A,B_{2})_{1}$ is $F-\delta$ multiplicative.\\
(3) For $g\in\widetilde{E}^{i}$ with $Y_{i}=T_{II,k}, T_{III,k}$ or $S^{2}$,
$$\|(\beta\oplus r)(\alpha_1(g))-(\psi\oplus(r\circ\alpha_{1}))(g)\|<5\varepsilon_{1};$$
for $g\in E^{i}(\supset\widetilde{E}^{i})$ with $Y_{i}=\{pt\}, [0,1], S^{1}$,
$$\|(\beta\oplus r)(\alpha_1(g))-(\psi\oplus r\circ\alpha_{1})(g)\|<\varepsilon_{1};$$
for all $f\in F$,
$$\|(\alpha_{2}\oplus id)\circ(\beta\oplus r)(f)-(\phi\oplus r)(f)\|<\varepsilon_{1}.$$
(4) $\omega(\beta(F)\cup\psi(E))<\varepsilon_{2}$.

The following lemma will be used in the proof of our main theorem.

\noindent\textbf{Lemma 3.15}~~Let $lim(A_{n}=\bigoplus^{k_{n}}_{i=1}M_{[n,i]}C(X_{n,i}),\phi_{n,m})$
be an $AH$ inductive system for which the  limit $C^{*}$-algebra has the ideal property, with $X_{n,i}$ connected simplicial
complexes and with $\sup\limits_{n,i}dim(X_{n,i})<+\infty$.
Let $P_{0}\in A_{n}$ be a trivial projection. Let $A=P_{0}A_{n}P_{0}(=\bigoplus_{i}P^{i}_{0}A^{i}_{n}P^{i}_{0})$.
Let $A^{'}=M_{L}(A), r(A)$ be as in 3.13 and 3.14 with $r:A\longrightarrow r(A)$. Let $\phi:A\longrightarrow A^{'}$ be a unital simple embedding (see 2.23).
Then there exist $m>n$, and a unital homomorphism
$$\lambda:M_{L}(A)\oplus r(A)\longrightarrow\phi_{n,m}(P_{0})A_{m}\phi_{n,m}(P_{0})$$
such that $\lambda\circ(\phi\oplus r)$ is homotopic to $\phi_{n,m}|_{P_{0}A_{n}P_{0}}:$
$P_{0}A_{n}P_{0}\longrightarrow\phi_{n,m}(P_{0})A_{m}\phi_{n,m}(P_{0})$.
That is the following diagram commutes up to homotopy
$$
\xymatrix{
A(=P_{0}A_{n}P_{0})\ar[rr]^{\phi_{n,m}}\ar[rdrd]^{\phi\oplus r}  &  &  \phi_{n,m}(P_{0})A_{m}\phi_{n,m}(P_{0})  &  &  \\
     &  &  &  &  &  &  \\
                                                                     &  &   M_{L}(A)\oplus r(A)\ar[uu]_{\lambda} & &   \\
 }
$$
\begin{proof}
One can construct the maps inside each block for each $\phi^{i,j}_{n,m}:P^{i}_{0}A^{i}_{n}P^{i}_{0}\longrightarrow\phi^{i,j}_{n,m}(P^{i}_{0})A^{j}_{m}\phi^{i,j}_{n,m}(P^i_0)$.
Then the lemma follows from Lemma 1.6.30 of [G4] and Lemma 2.20. Note that we use the fact that
if $\xi:M_{\bullet}(C(X))\longrightarrow B$ has finite dimensional image, then there is a map $\xi^{'}:r(M_{\bullet}C(X))\longrightarrow B$ such that $\xi^{'}\circ r$ is homotopic to $\xi$.\\
\end{proof}
\noindent\textbf{\S4. The proof of the main theorem}

The following lemma is one of the versions of the Elliott intertwining argument (see [Ell1]) or Proposition 3.1 of [D].

\noindent\textbf{Proposition 4.1}~~Consider the diagram
$$
\xymatrix{
A_{1}\ar[rr]^{\phi_{1,2}}\ar[rdrd]^{\beta_{1}}  &  &  A_{2}\ar[rr]^{\phi_{2,3}}\ar[rdrd]^{\beta_{2}}  &  &  A_{3}\ar[rr]
  &  & \cdot\cdot\cdot\ar[rr] &  & A_{n}\ar[rr]^{\phi_{n,n+1}}\ar[rdrd]^{\beta_{n}}   &  &   \cdot\cdot\cdot &  & \\
     &  &  &  &  &  &  \\
     B_{1}\ar[rr]_{\psi_{1,2}}\ar[uu]_{\alpha_{1}} &  &  B_{2}\ar[rr]_{\psi_{2,3}}
     \ar[uu]^{\alpha_{2}}  &  & B_{3} \ar[rr] &  & \cdot\cdot\cdot\ar[rr] &  & B_{n}\ar[rr]_{\psi_{n,n+1}}\ar[uu]^{\alpha_{n}}   &  & \ar[uu]_{\alpha_{n+1}}\cdot\cdot\cdot &  & \\
 }
$$
where $A_{n}$,$B_{n}$ are $C^{*}$-algebras $\phi_{n,n+1}, \psi_{n,n+1}$ are homomorphisms and $\alpha_{n}, \beta_{n}$ are linear $*$-contraction. Suppose that $F_{n}\subset A_{n}, \widetilde{E_{n}}\subset E_{n}\subset B_{n}$ are finite sets satisfying the following condition
$$\phi_{n,n+1}(F_{n})\cup\alpha_{n+1}(E_{n+1})\subset F_{n+1}, \psi_{n,n+1}(E_{n})\cup\beta_{n}(F_{n})\subset \widetilde{E_{n+1}}$$
 and
$\overline{\bigcup^{\infty}_{n=1}\phi_{n,\infty}(F_{n})}$
and
$\overline{\bigcup^{\infty}_{n=1}\psi_{n,\infty}(E_{n})}$ are the unit balls of $A=lim(A_{n},\phi_{n,m})$ and $B=lim(B_{n},\psi_{n,m})$, respectively. Suppose that there is a sequence $\varepsilon_{1}, \varepsilon_{2}...$ of positive numbers with $\sum\varepsilon_{n}<+\infty$ such that $\alpha_{n}$ and $\beta_{n}$ are $E_{n}-\varepsilon_{n}$ and $F_{n}-\varepsilon_{n}$ multiplicative, respectively, and $$\|\phi_{n,n+1}(f)-\alpha_{n+1}\circ\beta_{n}(f)\|<\varepsilon_{n},$$
and
$$\|\psi_{n,n+1}(g)-\beta_{n}\circ\alpha_{n}(g)\|<\varepsilon_{n},$$
for all $f\in F_{n}$ and $g\in\widetilde{E_{n}}$.\\
Then $A$ is isomorphic to $B$.\\

The following is the main theorem of this article.

\noindent\textbf{Theorem 4.2}~~Suppose that $lim(A_{n}=\bigoplus^{t_{n}}_{i=1}M_{[n,i]}C(X_{n,i}),\phi_{n,m})$ is an $AH$ inductive limit with $dim(X_{n,i})\leq M$ for a fixed positive integer $M$ such that the limit algebra has the ideal property. Then there is another inductive limit system$(B_{n}=\bigoplus^{s_{n}}_{i=1}M_{\{n,i\}}C(Y_{n,i}), \psi_{n,m})$ with the same limit algebra as the above system, where all the $Y_{n,i}$ are spaces of the form $\{pt\}, [0,1], S^{1}, S^{2}$, $T_{\uppercase\expandafter{\romannumeral2},k}$, $T_{\uppercase\expandafter{\romannumeral3},k}$.
 \begin{proof}
 Without loss of generality, assume that the spaces $X_{n,i}$ are connected finite simplicial complexes and $\phi_{n,m}$ are injective (see [EGL2]).

Let $\varepsilon_{1}>\varepsilon_{2}>\varepsilon_{3}>\cdot\cdot\cdot>0$ be a sequence of positive numbers satisfying $\sum\varepsilon_{n}<+\infty$. We need to construct the intertwining diagram
 $$
 \xymatrix@R=0.5ex{
 F_{1} & & F_{2} & & & & F_{n} & & F_{n+1}\\
 \bigcap & & \bigcap & & & & \bigcap & & \bigcap\\
 A_{s(1)} \ar[rdrd]^{\beta_{1}} \ar[rr]^{\phi_{s(1),s(2)}} & & A_{s(2)} \ar[rdrd]^{\beta_{2}} \ar[rr]^{\phi_{s(2),s(3)}} & &\cdots \ar[rr] & & A_{s(n)} \ar[rdrd]^{\beta_{n}} \ar[rr]^{\phi_{s(n),s(n+1)}}  & & A_{s(n+1)} \ar[rr] \ar[rdrd] & & \cdots\\
   & & & & & &\\
 B_{1} \ar[rr]^{\psi_{1,2}} \ar[uu]^{\alpha_{1}} & & B_{2} \ar[uu]^{\alpha_{2}} \ar[rr]^{\psi_{2,3}} & &\cdots \ar[rr] & & B_{n}\ar[rr]^{\psi_{n,n+1}} \ar[uu]^{\alpha_{n}} &  & B_{n+1} \ar[uu]^{\alpha_{n+1}} \ar[rr] & & \cdots\\
 \bigcup & & \bigcup & & & & \bigcup & & \bigcup\\
 E_{1} & & E_{2} & & & & E_{n} & & E_{n+1}\\
 \bigcup & & \bigcup & & & & \bigcup & & \bigcup\\
 \widetilde{E_{1}} & & \widetilde{E_{2}} & & & & \widetilde{E_{n}} & & \widetilde{E_{n+1}}
 }
 $$

satisfying the following conditions\\
(0.1) $(A_{s(n)},\phi_{s(n),s(m)})$ is a subinductive system of $(A_{n}, \phi_{n,m})$, and $(B_{n}, \psi_{n,m})$ is an inductive system of matrix algebras over the spaces $\{pt\}, [0,1], S^{1}, {T_{\uppercase\expandafter{\romannumeral2},k}} , {T_{\uppercase\expandafter{\romannumeral3},k}},S^{2}.$\\
(0.2) Choose $\{a_{ij}\}^{\infty}_{j=1}\subset A_{s(i)}$ and $\{b_{ij}\}^{\infty}_{j=1}\subset B_{i}$ to be countable dense subsets of the unit balls of $A_{s(i)}$ and $B_{i}$, respectively. $F_{n}$ are subsets of the unit balls of $A_{s(n)}$, and $\widetilde{E_{n}}\subset E_{n}$ are both subsets of the unit balls of $B_{n}$ satisfying
$$\phi_{s(n),s(n+1)}(F_{n})\cup\alpha_{n+1}(E_{n+1})\cup\bigcup^{k+1}_{i=1}\phi_{s(i),s(n+1)}(\{a_{i1},a_{i2}\cdot\cdot\cdot a_{in+1}\})\subset F_{n+1}$$
$$\psi_{n,n+1}(E_{n})\cup\beta_{n}(F_{n})\subset \widetilde{E}_{n+1}\subset E_{n+1}$$
and
$$\bigcup^{n+1}_{i=1}\psi_{i,n+1}(\{b_{i1},b_{i2}\cdot\cdot\cdot b_{in+1}\})\subset E_{n+1}$$
(0.3) $\beta_{n}$ are $F_{n}-2\varepsilon_{n}$ multiplicative and $\alpha_{n}$ are homomorphisms\\
(0.4) $\|\psi_{n,n+1}(g)-\beta_{n}\circ\alpha_{n}(g)\|<8\varepsilon_{n}$ for all $g\in\widetilde{E}_{n}$ and $\|\phi_{s(n),s(n+1)}(f)-\alpha_{n+1}\circ\beta_{n}(f)\|<14\varepsilon_{n}$ for all $f\in F_{n}$\\
(0.5) For any block $B^{i}_{n}$ with spectrum $T_{\uppercase\expandafter{\romannumeral2},k},T_{\uppercase\expandafter{\romannumeral3},k},S^{2},\omega(\widetilde{E}_{n}^{i})<\varepsilon_{n}$, where $\widetilde{E}_{n}^{i}=\pi_{i}(\widetilde{E}_{n})$ for $\pi_{i}:B_{n}\longrightarrow B^{i}_{n}$ the canonical projections.

The diagram will be constructed inductively. First, let $B_{1}=\{0\}, A_{s(1)}=A_{1}, \alpha_{1}=0$. Let $b_{1j}=0\in B_{1}$ for $j=1,2,...$ and let $\{a_{1j}\}^{\infty}_{j=1}$ be a countable dense subset of the unit ball of $A_{s(1)}$. And let $\widetilde{E}_{1}=E_{1}=\{b_{11}\}=B_{1}$
and $F_{1}=\bigoplus^{t_{1}}_{i=1}F^{i}_{1}$, where $F^{i}_{1}=\pi_{i}(\{a_{11}\})\subset A^{i}_{1}$.\\
As inductive assumption, assume that we already have the diagram

 $$
 \xymatrix@R=0.5ex{
F_{1} & & F_{2} & & & & F_{n}\\
\bigcap & & \bigcap & & & & \bigcap\\
 A_{s(1)} \ar[rdrd]^{\beta_{1}} \ar[rr]^{\phi_{s(1),s(2)}} & & A_{s(2)} \ar[rdrd]^{\beta_{2}} \ar[rr]^{\phi_{s(2),s(3)}} & &\cdots \ar[rr] \ar[rdrd]^{\beta_{n-1}} & & A_{s(n)}\\
   & & & & & &\\
 B_{1} \ar[rr]^{\psi_{1,2}} \ar[uu]^{\alpha_{1}} & & B_{2} \ar[uu]^{\alpha_{2}} \ar[rr]^{\psi_{2,3}} & &\cdots \ar[rr] & & B_{n}\ar[uu]^{\alpha_{n}}\\
 \bigcup & & \bigcup & & & & \bigcup\\
 E_{1} & & E_{2} & & & & E_{n}\\
 \bigcup & & \bigcup & & & & \bigcup\\
 \widetilde{E_{1}} & & \widetilde{E_{2}} & & & & \widetilde{E_{n}}
 }
 $$
and for each $i=1,2,\cdot\cdot\cdot,n$, we have dense subsets
$\{a_{ij}\}^{\infty}_{j=1}\subset$ the unit ball of $A_{s(i)}$
and
$\{b_{ij}\}^{\infty}_{j=1}\subset$ the unit ball of $B_{i}$
 satisfying the conditions (0.1)-(0.5) above. We have to construct the next piece of the diagram
$$
\xymatrix@!C{
F_{n}\subset A_{s(n)}\ar[rr]^{\phi_{s(n),s(n+1)}}\ar[rdrd]^{\beta_{n}}  &  &  A_{s(n+1)}\supset F_{n+1}  &  &    \\
     &  &  &  &  &  &  \\
\widetilde{E}_{n}\subset E_{n}
\subset B_{n}\ar@<-8mm>[uu]^{\alpha_{n}}\ar[rr]_{\psi_{n,n+1}} &  &
   \;\;B_{n+1}\ar@<10mm>[uu]_{\alpha_{n+1}}\supset
E_{n+1}\supset\widetilde{E}_{n+1}  &  &  \\
 }
$$
to satisfy the conditions (0.1)-(0.5).

Among the conditions for the  induction assumption, we will only use the conditions that $\alpha_{n}$ is a homomorphism and (0.5) above.\\
\textbf{Step 1}. We enlarge $\widetilde{E}_{n}$ to $\bigoplus_{i}\pi_{i}(\widetilde{E}_{n})$ and $E_{n}$ to $\bigoplus_{i}\pi_{i}(E_{n})$. Then
$$\widetilde{E}_{n}(=\bigoplus\widetilde{E}_{n}^{i})\subset E_{n}(=\bigoplus E_{n}^i)$$ and for each $B^{i}_{n}$ with spectrum $T_{\uppercase\expandafter{\romannumeral2},k},T_{\uppercase\expandafter{\romannumeral3},k},S^{2}$, we have $\omega(\widetilde{E}_{n}^{i})<\varepsilon_{n}$ from the induction assumption (0.5). By Theorem 3.9 applied to $\alpha_{n}:B_{n}\rightarrow A_{s(n)},\widetilde{E}_{n}\subset E_{n}\subset B_{n},F_{n}\subset A_{s(n)}$ and $\varepsilon_{n}>0$, there are $A_{m_{1}}(m_{1}>s(n))$, two orthogonal projections $P_{0},P_{1}\in A_{m_{1}}$ with $\phi_{s(n),m_{1}}(\textbf{1}_{A_{s(n)}})=P_{0}+P_{1}$ and $P_{0}$ trivial, a $C^{*}$-algebra $C$---a direct sum of matrix algebras over $C[0,1]$ or $\mathbb{C}$, a unital map $\theta\in Map(A_{s(n)},P_{0}A_{m_{1}}P_{0})_{1}$, a unital homomorphism $\xi_{1}\in Hom(A_{s(n)},C)_{1}$, a unital homomorphism $\xi_{2}\in Hom(C,P_{1}A_{m_{1}}P_{1})_{1}$ and a homomorphism $\alpha\in Hom(B_n, P_0A_{m_1}P_0)$ such that\\
(1.1) $\|\phi_{s(n),m_{1}}(f)-\theta(f)\oplus(\xi_{2}\circ\xi_{1})(f)\|<\varepsilon_{n}$ for all $f\in F_{n}$.\\
(1.2) $\theta$ is $F_{n}-\varepsilon_{n}$ multiplicative and $F:=\theta(F_{n})$ satisfies  $\omega(F)<\varepsilon_{n}$.\\
(1.3) $\|\alpha(g)-\theta\circ\alpha_{n}(g)\|<3\varepsilon_{n}$ for all $g\in\widetilde{E}_{n}$.

Let all the blocks of $C$ be parts of  the $C^{\ast}$-algebra $B_{n+1}$. That is
$$B_{n+1}=C\oplus(some\; other\; blocks).$$
The map $\beta_{n}:A_{s(n)}\rightarrow B_{n+1}$, and the homomorphism $\psi_{n,n+1}:B_{n}\rightarrow B_{n+1}$ are defined by $\beta_{n}=\xi_{1}:A_{s(n)}\rightarrow C(\subset B_{n+1})$
and $\psi_{n,n+1}=\xi_{1}\circ\alpha_{n}:B_{n}\rightarrow C(\subset B_{n+1})$
for the blocks of $C(\subset B_{n+1})$. For this part, $\beta_{n}$ is also a homomorphism.\\
\textbf{Step 2}. Let $A=P_{0}A_{m_{1}}P_{0},F=\theta(F_{n})$. Since $P_{0}$ is a trivial projection, $$A\cong\oplus M_{l_{i}}(C(X_{m_{1},i})).$$ Let $r(A):=\bigoplus M_{l_{i}}(\mathbb{C})$ and $r:A\rightarrow r(A)$ be as in 3.13. Applying Corollary 3.14 to $\alpha:B_{n}\rightarrow A,\widetilde{E}_{n}\subset E_{n}\subset B_{n}$ and $F\subset A$, we obtain the following diagram
$$
\xymatrix{
A\ar[rr]^{\phi\oplus r}\ar[rdrd]^{\beta}      & &M_{L}(A)\oplus r(A)  &  &    \\
     &  &  &  &  &  &  \\
B_{n}\ar[rr]^{\psi}\ar[uu]^{\alpha} &  &  B\ar[uu]_{\alpha{'}} &  &  \\
 }
$$
such that\\
(2.1) $B$ is a direct sum of matrix algebras over $\{pt\},[0,1],S^{1},T_{\uppercase\expandafter{\romannumeral2},k},T_{\uppercase\expandafter{\romannumeral3},k}$ or $S^{2}$.\\
(2.2) $\alpha^{\prime}$ is an injective homomorphism and $\beta$ is $F-\varepsilon_{n}$ multiplicative.\\
(2.3) $\phi:A\rightarrow M_{L}(A)$ is a unital simple embedding and $r:A\rightarrow r(A)$ is as in 3.13.\\
(2.4) $\|\beta\circ\alpha(g)-\psi(g)\|<5\varepsilon_{n}$ for all $g\in \widetilde{E}_{n}$ and $\|(\phi\oplus r)(f)-\alpha{'}\circ\beta(f)\|<\varepsilon_{n}$ for all $f\in F(:=\theta(F_{n}))$.\\
(2.5) $\omega(\psi(E_{n})\cup\beta(F))<\varepsilon_{n+1}$ (note that $\beta(F)=\beta\circ\theta(F_{n})$).

Let all the blocks B be also part of $B_{n+1}$, that is
$$B_{n+1}=C\oplus B\oplus(some\; other\; blocks)$$
The maps $\beta_{n}:A_{s(n)}\longrightarrow B_{n+1},\psi_{n,n+1}:B_{n}\longrightarrow B_{n+1}$
are defined by
$$\beta_{n}:=\beta\circ\theta:A_{s(n)}\xrightarrow{\theta}A\xrightarrow{\beta}B(\subset B_{n+1})$$
and
$$\psi_{n, n+1}:=\psi:B_{n}\rightarrow B(\subset B_{n+1})$$ for the blocks of $B(\subset B_{n+1})$. This part of $\beta_{n}$ is $F_{n}-2\varepsilon_{n}$ multiplicative, since $\theta$ is $F_{n}-\varepsilon_{n}$ multiplicative, $\beta$ is $F-\varepsilon_{n}$ multiplicative and $F=\theta(F_{n})$.\\
\textbf{Step 3}. By Lemma 3.15 applied to $\phi\oplus r:A\rightarrow M_{L}(A)\oplus r(A)$, there is an $A_{m_{2}}$ and there is a unital homomorphism
$$\lambda:M_{L}(A)\oplus r(A)\rightarrow RA_{m_{2}}R,$$ where $R=\phi_{m_{1},m_{2}}(P_{0})$ (write $R$ as $\bigoplus_{j}R^{j}\in\bigoplus_{j}A^{j}_{m}$) such that the diagram
$$
\xymatrix{
A(=P_{0}A_{m}P_{0})\ar[rr]^{\phi_{m_{1},m_{2}}}\ar[rdrd]^{\phi\oplus r}  &  &  RA_{m_{2}}R  &  &  \\
     &  &  &  &  &  &  \\
                                                                     &  &   M_{L}(A)\oplus r(A)\ar[uu]_{\lambda} & &   \\
 }
$$
satisfies the following condition:\\
(3.1) $\lambda\circ(\phi\oplus r)$ is homotopy equivalent to
$$\phi{'}:=\phi_{m_{1},m_{2}}|_{A}.$$
\textbf{Step 4}. Applying Theorem 1.6.9 of [G4] to the finite set $F\subset A$ (with $\omega(F)<\varepsilon_{n}$)
and to the two homotopic homomorphisms $\phi{'}$ and $\lambda\circ(\phi\oplus r):A\rightarrow RA_{m_{2}}R$ (with $RA_{m_{2}}R$ in place of $C$ in Theorem 1.6.9 of [G4]), we obtain a finite set $F{'}\subset RA_{m_{2}}R, \delta>0$ and $L>0$ as in the Theorem.

Let $G:=\psi(E_{n})\cup\beta(F)$. From (2.5), we have $\omega(G)<\varepsilon_{n+1}<\varepsilon_{n}$. By Theorem 3.8 applied to $RA_{m_{2}}R$ and
$$\lambda\circ\alpha{'}:B\rightarrow RA_{m_{2}}R$$ finite set $G\subset B$, $F{'}\cup\phi^{\prime}(F)\subset RA_{m_{2}}R$, $min \{\varepsilon_{n},\delta\}>0$ (in place of $\varepsilon$) and $L>0$, there are $A_{s(n+1)}$, mutually orthogonal projections $Q_{0},Q_{1},Q_{2}\in A_{s(n+1)}$ with $\phi_{m_{2},s(n+1)}(R)=Q_{0}\oplus Q_{1}\oplus Q_{2}$, a $C^{*}$-algebra $D$---a direct sum of matrix algebras over C[0,1] or $\mathbb{C}$---,
a unital map $\theta_{0}\in$ Map$(RA_{m_{2}}R,Q_{0}A_{s(n+1)}Q_{0})_1$, and four unital homomorphisms $$\theta_{1}\in Hom(RA_{m_{2}}R,Q_{1}A_{s(n+1)}Q_{1})_{1},\xi_{3}\in Hom(RA_{m_{2}}R,D)_{1},\xi_{4}\in Hom(D,Q_{2}A_{s(n+1)}Q_{2})_{1}$$ and $\alpha{''}\in Hom(B,(Q_{0}+Q_{1})A_{s(n+1)}(Q_{0}+Q_{1}))_{1}$ such that the following statements are true.\\
(4.1) $\|\phi_{m_{2},s(n+1)}(f)-((\theta_{0}+\theta_{1})\oplus\xi_{4}\circ\xi_{3})(f)\|<\varepsilon_{n}$ for all $f\in\phi_{m_{1},m_{2}}|_{A}(F)\subset RA_{m_{2}}R$.\\
(4.2) $\|\alpha{''}(g)-(\theta_{0}+\theta_{1})\circ\lambda\circ\alpha{'}(g)\|<3\varepsilon_{n+1}<3\varepsilon_{n}$
for all $g\in G$.\\
(4.3) $\theta_{0}$ is $F{'}-min(\varepsilon_{n},\delta)$ multiplicative and $\theta_{1}$ satisfies that
$$\theta^{i,j}_{1}([q])>L\cdot[\theta^{i,j}_{0}(R^{i})]$$
for any non zero projection
$q\in R^{i}A_{m_{1}}R^{i}$.\\
By Theorem 1.6.9 of [G4], there is a unitary $u\in(Q_{0}\oplus Q_{1})A_{s(n+1)}(Q_{0}\oplus Q_{1})$ such that $$\|(\theta_{0}+\theta_{1})\circ\phi{'}(f)-Adu\circ(\theta_{0}+\theta_{1})\circ\lambda\circ(\phi\oplus r)(f)\|<8\varepsilon_{n}$$
for all $f\in F$.\\
Combining with second inequality of (2.4), we have\\
(4.4) $\|(\theta_{0}+\theta_{1})\circ\phi{'}(f)-Adu\circ(\theta_{0}+\theta_{1})\circ\lambda\circ\alpha{'}\circ\beta(f)\|<9\varepsilon_{n}$
for all $f\in F$.\\
\textbf{Step 5}. Finally let all blocks of $D$ be the rest of $B_{n+1}$. Namely, let $$B_{n+1}=C\oplus B\oplus D,$$ where $C$ is from Step 1, $B$ is from Step 2 and $D$ is from Step 4.

We already have the definition of $\beta_{n}:A_{s(n)}\rightarrow B_{n+1}$ and $\psi_{n,n+1}: B_{n}\rightarrow B_{n+1}$
for those blocks of $C\oplus B\subset B_{n+1}$
(from Step 1 and Step 2). The definition of $\beta_{n}$ and $\psi_{n,n+1}$ for blocks of $D$ and the homomorphism $\alpha_{n+1}:C\oplus B\oplus D\rightarrow A_{s(n+1)}$
will be given below.

The part of $\beta_{n}:A_{s(n)}\rightarrow D(\subset B_{n+1})$
is defined by
$$\beta_{n}=\xi_{3}\circ\phi^{'}\circ\theta:A_{s(n)}\xlongrightarrow{\theta}A\xlongrightarrow{\phi^{\prime}}RA_{m_{2}}R\xlongrightarrow{\xi_{3}}D$$
(Recall that $A=P_{0}A_{m_{2}}P_{0}$ and $\phi{'}=\phi_{m_{1},m_{2}}|_{A})$. Since $\theta$ is $F_{n}-\varepsilon_{n}$ multiplicative, and $\phi{'}$ and $\xi_{3}$ are homomorphisms, we know this part of $\beta_{n}$ is $F_{n}-\varepsilon_{n}$ multiplicative.

The part of $\psi_{n,n+1}:B_{n}\rightarrow D(\subset B_{n+1})$ is defined by
$$\psi_{n,n+1}=\xi_{3}\circ\phi{'}\circ\alpha:B_{n}\xlongrightarrow{\alpha}A\xlongrightarrow{\phi{'}}RA_{m}R\xlongrightarrow{\xi_{3}}D$$
which is a homomorphism.

The homomorphism $\alpha_{n+1}:C\oplus B\oplus D\rightarrow A_{s(n+1)}$
is defined as following.

Let $\phi{''}=\phi_{m_{1},s(n+1)}|_{P_{1}A_{m_{1}}P_{1}}:P_{1}A_{m_{1}}P_{1}\longrightarrow
\phi_{m_{1},s(n+1)}(P_{1})A_{s(n+1)}\phi_{m_{1},s(n+1)}(P_{1})$, where $P_{1}$ is from Step 1.
Define $$\alpha_{n+1}|_{C}=\phi{''}\circ\xi_{2}:C\xlongrightarrow{\xi_{2}}P_{1}A_{m_{1}}P_{1}\xlongrightarrow{\phi{''}}\phi_{m_{1},s(n+1)}(P_{1})A_{s(m+1)}\phi_{m_{1},s(n+1)}(P_{1})$$
where $\xi_{2}$ is from Step 1.\\
\begin{center}
$\alpha_{n+1}|_{B}=Adu\circ\alpha{''}:B\xlongrightarrow{\alpha{''}}(Q_{0}\oplus Q_{1})A_{s(n+1)}(Q_{0}+Q_{1})\xlongrightarrow{Adu}(Q_{0}\oplus Q_{1})A_{s(n+1)}(Q_{0}+Q_{1})$,
\end{center}
where $\alpha{''}$ is from Step 4, and define
$$\alpha_{n+1}|_{D}=\xi_{4}:D\rightarrow Q_{2}A_{s(n+1)}Q_{2}.$$
Finally choose $\{a_{n+1,j}\}^{\infty}_{j=1}\subset A_{s(n+1)}$ and $\{b_{n+1,j}\}^{\infty}_{j=1}\subset B_{n+1}$
to be countable dense subsets of the unit balls of $A_{s(n+1)}$ and $B_{n+1}$, respectively. And choose\\
$$F{'}_{n+1}=\phi_{s(n),s(n+1)}(F_{n})\cup\alpha_{n+1}(E_{n+1})\cup\bigcup^{n+1}_{i=1}\phi_{s(i),s(n+1)}(\{a_{i1},a_{i2},\cdot\cdot\cdot, a_{in+1}\})$$\\
$$E{'}_{n+1}=\psi_{n,n+1}(E_{n})\cup\beta_{n}(F_{n})\cup\bigcup^{n+1}_{i=1}\psi_{i,n+1}(\{b_{i1},b_{i2},\cdot\cdot\cdot, b_{in+1}\})$$\\
$$\widetilde{E}'_{n+1}=\psi_{n,n+1}(E_{n})\cup\beta_{n}(F_{n})\subset E'_{n+1}.$$
Define $F^{i}_{n+1}=\pi_{i}(F^{\prime}_{n+1})$ and $F_{n+1}=\bigoplus_{i}F^{i}_{n+1}$, $E^{i}_{n+1}=\pi_{i}(E{'}_{n+1})$
and $E_{n+1}=\oplus_{i}E^{i}_{n+1}$.
For those blocks $B^{i}_{n+1}$ inside the algebra $B$ define $\widetilde{E}_{n+1}^{i}=\pi_{i}(\widetilde{E}'_{n+1})$.
For those blocks inside $C$ and $D$, define $\widetilde{E}_{n+1}^{i}=E^{i}_{n+1}$.
And finally let $\widetilde{E}_{n+1}=\bigoplus_{i}\widetilde{E}_{n+1}^{i}$.
Note that all the blocks with spectrum $T_{\uppercase\expandafter{\romannumeral2},k},T_{\uppercase\expandafter{\romannumeral3},k}$ and $S^{2}$ are in $B$. And hence (2.5) tells us that for those blocks, $\omega(\widetilde{E}_{n+1}^{i})<\varepsilon_{n+1}$.\\
Thus we obtain the following diagram
$$
\xymatrix@!C{
F_{n}\subset A_{s(n)}\ar[rr]^{\phi_{s(n),s(n+1)}}\ar[drdr]^{\beta_{n}}  &  &  A_{s(n+1)}\supset F_{n+1}  &  &    \\
     &  &  &  &  &  &  \\
\widetilde{E}_{n}\subset E_{n}
\subset B_{n}\ar@<-8mm>[uu]^{\alpha_{n}}\ar[rr]_{\psi_{n,n+1}} &  &
   \;\;B_{n+1}\ar@<12mm>[uu]_{\alpha_{n+1}}\supset
E_{n+1}\supset\widetilde{E}_{n+1}.  &  &  \\
 }
$$
\textbf{Step 6}.~~Now we need to verify all the conditions (0.1)-(0.5) for the above diagram.

From the end of Step 5, we know (0.5) holds, (0.1)-(0.2) hold from the construction (see the construction of $B, C, D$ in Step 1, 2 and 4, and $\widetilde{E}_{n+1}\subset E_{n+1}, F_{n+1}$ is the end of Step 5).
(0.3) follows from the end of Step 1, the end of Step 2 and the part of definition of $\beta_{n}$ for $D$ from Step 5.

So we only need to verify (0.4).

Combining (1.1) with (4.1), we have
$$\|\phi_{s(n),s(n+1)}(f)-[(\phi^{''}\circ\xi_{2}\circ\xi_{1})\oplus(\theta_{0}+\theta_{1})\circ\phi^{'}\circ\theta\oplus(\xi_{4}\circ\xi_{3}\circ\phi^{'}\circ\theta)](f)\|
<\varepsilon_{n}+\varepsilon_{n}=2\varepsilon_{n}$$
for all $f\in F_{n}$ (recall that $\phi^{''}=\phi_{m_{1},s(n+1)}|_{P_{1}A_{m_{1}}P_{1}},\phi^{'}:=\phi_{m_{1},m_{2}}|_{P_{0}A_{m_{1}}P_{0}}$).

Combining with (4.2) and (4.4), and the definitions of $\beta_{n}$ and $\alpha_{n+1}$, the above inequality yields
$$\|\phi_{s(n),s(n+1)}(f)-(\alpha_{n+1}\circ\beta_{n})(f)\|< 9\varepsilon_{n}+3\varepsilon_{n}+2\varepsilon_{n}=14\varepsilon_{n},~~~\forall f\in F_{n}.$$
Combining (1.3), the first inequality of (2.4) and the definition of $\beta_{n}$ and $\psi_{n,n+1}$, we have
$$\|\psi_{n,n+1}(g)-(\beta_{n}\circ\alpha_{n})(g)\|<5\varepsilon_{n}+3\varepsilon_{n}=8\varepsilon_{n},~~~\forall g\in \widetilde{E}_{n}.$$
So we obtain (0.4).\\
The theorem follows from Proposition 4.1.\\
\end{proof} 

\vspace{0.3in}


Guihua Gong, College of Mathematics and Information Science, Hebei Normal University,  Shijiazhuang, Hebei, 050024,    China, and\\
Department of Mathematics, University of Puerto Rico at Rio Piedras, PR 00936, USA\\
email address:  guihua.gong@upr.edu\\

Chunlan Jiang, College of Mathematics and Information Science, Hebei Normal University,  Shijiazhuang, Hebei, 050024,    China\\
email address:  cljiang@hebtu.edu.cn\\

Liangqing Li, Department of Mathematics, University of Puerto Rico at Rio Piedras, PR 00936, USA\\
email address:  liangqing.li@upr.edu\\

Cornel Pasnicu, Department of Mathematics, University of Texas at San Antonio, San Antonio,  TX 78249, USA\\
email address:  Cornel.Pasnicu@utsa.edu

\clearpage

\begin{tiny}

\begin{small}

\end{small}

\end{tiny}

\end{document}